\documentclass[3p]{elsarticle}
\usepackage[utf8]{inputenc}
\usepackage[ngerman, english]{babel}
\usepackage{amsfonts,amsmath,amssymb,amsthm,mathtools,esint}
\usepackage{algorithm}
\usepackage[hidelinks]{hyperref}
\usepackage{cleveref}
\usepackage{algpseudocode}
\usepackage{placeins}
\usepackage{fixme}
\usepackage{graphicx}
\usepackage{caption}
\usepackage[utf8]{inputenc}
\usepackage[T1]{fontenc}
\usepackage{graphicx}
\usepackage{amsmath}
\usepackage{amssymb}
\usepackage{amsthm}
\usepackage{mathtools}
\usepackage{listings}
\usepackage{enumerate}
\usepackage[autostyle=true,german=quotes]{csquotes}
\usepackage{setspace}

\usepackage{paralist}
\usepackage{bbm}
\usepackage{color}
\usepackage{float}

\usepackage{tikz}
\usepackage{pgfplots,xparse}

\usetikzlibrary{calc,shapes,patterns,spy,external}
\tikzset{external/optimize=false}
\tikzset{external/system call= {pdflatex -save-size=80000 
                           -pool-size=10000000 
                           \tikzexternalcheckshellescape 
                           -halt-on-error 
                           -interaction=batchmode
                           -jobname "\image" "\texsource"}} 
\makeatletter
\newcommand{\useexternalfile}[1]{%
  \tikzsetnextfilename{#1_out}%
  \input{\tikzexternal@filenameprefix#1.tex}}
\makeatother

\tikzset{slopetriangle/.style={%
  bottom color=black!20,
  middle color=black!5,
  top color=white,
  draw=black
}}

\pgfplotsset{%
  width=.65\linewidth,
  axis background/.style={fill=black!5!white},
  grid style={densely dotted,semithick},
  legend style={%
    legend columns=1,
    legend pos=outer north east
  },
  compat=newest % compatibility for old pgfplots versions
}
\ExplSyntaxOn
\DeclareExpandableDocumentCommand \eval { m } { \fp_eval:n { #1 } }
\ExplSyntaxOff

\usetikzlibrary{external}
\definecolor{mygray}{rgb}{0.5,0.5,0.5}
\definecolor{mymauve}{rgb}{0.58,0,0.82}
\definecolor{mygreen}{RGB}{28,172,0} % color values Red, Green, Blue
\definecolor{mylilas}{RGB}{170,55,241}

\newcommand{\R}{{\mathbb R}}
\newcommand{\N}{{\mathbb N}}

\newcommand{\F}{\mathcal{F}}
\newcommand{\E}{\mathcal{E}}
\newcommand{\T}{\mathcal{T}}
\newcommand{\V}{\mathcal{V}}

\theoremstyle{remark}

\newtheorem*{rem}{Remark}

\theoremstyle{definition}

\theoremstyle{plain}
\newtheorem{theorem}{Theorem}[section]
\newtheorem{lemma}[theorem]{Lemma}

\newtheorem{prop}[theorem]{Proposition}
\newtheorem{cor}[theorem]{\protect\colname}

\newcommand{\colname}{}% initialization
\addto\captionsenglish{%
  \renewcommand{\colname}{Corollary}%
}
\addto\captionsngerman{%
  \renewcommand{\colname}{Korollar}%
}

\usepackage{mathrsfs}

\fxsetup{%
	draft,%
	layout = inline,%
	theme = color,
	silent = true,%
	mode=multiuser,
	}
\FXRegisterAuthor{bg}{envbg}{\color{blue} ToDo}

\usetikzlibrary{shapes, arrows}
\usetikzlibrary{positioning,backgrounds}
\usetikzlibrary{arrows}
\usetikzlibrary{decorations.markings}
\usetikzlibrary{patterns}
\usetikzlibrary{fit}
\newcommand{\Vertices}{\mathcal V}

\newcommand{\n}{\nu}
\newcommand{\nF}{\n_E}
\renewcommand{\t}{\tau}
\newcommand{\BiL}{\Delta^2}
\newcommand{\JF}[1]{[#1]_E}
\newcommand{\Dn}[1]{{\Dnb #1}}
\newcommand{\DnF}[1]{\partial_{\nF} #1}
\newcommand{\DnnF}[1]{\partial_{\n_E\n_E}^2 #1}
\newcommand{\DLn}[1]{\partial_{\n} \Delta #1}
\newcommand{\Dt}[1]{\Dtb #1}
\newcommand{\Dtb}[1]{\partial_{\t} #1}

\newcommand{\Dnb}[1]{\partial_{\n} #1}
\newcommand{\Dnn}[1]{\partial_{\n\n}^2 #1}

\newcommand{\Dtn}[1]{\partial_{\t\n}^2 #1}
\newcommand{\Dtt}[1]{\partial_{\t\t}^2 #1}
\newcommand{\Dttn}[1]{\partial_{\t\t\n}^3 #1}
\newcommand{\Dttt}[1]{\partial_{\t\t\t}^3 #1}
\newcommand{\Dtttt}[1]{\partial_{\t\t\t\t}^4 #1}

\newcommand{\LL}[1]{\Delta^2 #1}

\newcommand{\LFace}{\LEdge}
\newcommand{\LEdge}{{L^2(E)}}
\newcommand{\HtwoO}{{H^2(\Omega)}}
\newcommand{\LO}{{L^2(\Omega)}}

\newcommand{\GCS}{{\Gamma_C\cup\Gamma_S}}
\newcommand{\GF}{{\Gamma_F}}
\newcommand{\GS}{{\Gamma_S}}
\newcommand{\GC}{{\Gamma_C}}
\renewcommand{\V}{{V}}

\newcommand{\TT}{\mathbb T}
\newcommand{\Tt}{\hat \T}
\newcommand{\Ft}{\hat \F}

\newcommand{\VT}{\V(\T)}
\newcommand{\VTt}{\V(\Tt)}

\newcommand{\uh}{u_h}
\newcommand{\tuh}{\tilde u_h}
\newcommand{\ul}{u_\ell}
\newcommand{\vh}{v_h}
\newcommand{\uht}{\hat \uh}
\newcommand{\vht}{\hat \vh}

\newcommand{\capx}{c_{\rm apx}}

\newcommand{\cb}{c_{\rm b}}
\newcommand{\cc}{c_{\rm c}}

\newcommand{\Cbiharm}{C_{\rm osc}}

\newcommand{\It}{\hat\I}

\newcommand{\Argyris}{\mathfrak A}

\newcommand{\Astd}{\Argyris_{\rm std}}
\newcommand{\Aext}{\Argyris_{\rm ext}}
\newcommand{\Tl}{\T_{\ell}}
\newcommand{\Tlnext}{\T_{\ell+1}}
\newcommand{\et}{\hat e}
\newcommand{\vt}{\hat v}

\newcommand{\wg}{{e_{\rm b}}}
\newcommand{\wgt}{\widehat \wg}
\newcommand{\w}{{e_0}}
\newcommand{\wt}{\widehat \w}
\newcommand{\vl}{v_\ell}

\newcommand{\osc}{\mathrm{osc}}

\newcommand{\midE}{\mathrm{mid}\; E}
\newcommand{\midF}{\midE}

\newcommand{\Nodes}{\mathcal N}

\newcommand{\ndof}{\ensuremath{N}}

\newcommand{\opA}{A}
\newcommand{\opAl}{\opA_\ell}

\newcommand{\opQ}{P^\top}
\newcommand{\opQl}{\opQ_\ell}

\newcommand{\opB}{B}
\newcommand{\opBl}{\opB_\ell}
\newcommand{\opS}{S}
\newcommand{\opSl}{\opS_\ell}
\newcommand{\opI}{I}

\newcommand{\opIl}{\opI}
\newcommand{\ellminus}{_{\ell-1}}
\newcommand{\opIBAl}{\ensuremath{\opIl-\opBl\opAl}}

\newcommand{\eg}{e.g.,\ }

\newcommand{\El}{\E_\ell}

\newcommand{\etaalg}{\eta_{\rm alg}}

\newcommand{\EO}{\E_\Omega}
\newcommand{\VerticesO}{\Vertices_\Omega}

\renewcommand{\F}{\E}
\newcommand{\widebar}{\overline}

\renewcommand{\hat}{\widehat}

\newcommand{\ula}{u_{\ell,0}}
\newcommand{\tula}{\tilde u_{\ell, 0}}

\newcommand{\txl}{\tilde x_{\ell}}
\newcommand{\Il}{\I_\ell}

\newcommand{\vflj}{\varphi^\ell_j}
\newcommand{\nuz}{\zeta_z}
\newcommand{\tauz}{\xi_z}
\newcommand{\dof}{L}
\newcommand{\Lzj}{\dof_{z,j}}
\newcommand{\bzj}{\varphi_{z,j}}
\newcommand{\bE}{b_E}
\newcommand{\dE}{d_E}

\newcommand{\I}{\mathcal{I}}

\newcommand{\trb}[1]{|\!|\!|#1|\!|\!|}

\newcommand{\vfzj}{\varphi_{z,j}}
\newcommand{\uhz}{u_{h,0}}

\hypersetup{pdftitle={Optimal multilevel adaptive FEM for the Argyris element}, pdfauthor={Benedikt Gräßle}}

\begin{document}
\title{Optimal multilevel adaptive FEM for the Argyris element}%
\author{Benedikt Gräßle}
\ead{graesslb@math.hu-berlin.de}
	\address{%
	Institut f\"ur Mathematik,
	Humboldt-Universit\"at zu Berlin,
	10117 Berlin, Germany}
\date{\today}

\begin{abstract}
	\noindent The main drawback for the application of the conforming Argyris FEM is the labourious implementation on the one hand and the low convergence rates on the other.
If no appropriate adaptive meshes are utilised, only the convergence rate caused by corner singularities
[Blum and Rannacher, 1980],
far below the
approximation order for smooth functions, can be achieved.
The fine
approximation with the Argyris FEM produces high-dimensional linear systems and for a long time an optimal preconditioned
scheme was not available for unstructured
grids.
This paper presents numerical benchmarks to confirm that the adaptive multilevel solver for the hierarchical Argyris FEM from
[Carstensen and Hu, 2021] 
is in fact
highly efficient and of linear time complexity.
Moreover, the very first display of optimal convergence rates in practically relevant benchmarks with corner
singularities and general boundary conditions leads to the rehabilitation of the Argyris finite element from the
computational perspective.
\end{abstract}
\maketitle

\section{Introduction}
\label{sec:Intro}
This paper discusses numerical aspects of an adaptive multilevel algorithm based on the Argyris finite element for the
biharmonic plate problem
with inhomogeneous mixed boundary conditions.
\subsection{Motivation}%
\label{sub:Motivation}
The conforming discretisation of fourth-order problems on unstructured domains with the finite element method (FEM) requires complicated $C^1$ elements
 like 
Hsieh-Clough-Tocher or Argyris elements \cite{ciarlet_finite_2002}.
Classical a priori analysis yields optimal rates of convergence 
for sufficiently smooth solutions only.
In practical applications however, singularities in the data or boundary of the domain lead to singular solutions
\cite{blum_boundary_1980}
and reduced convergence rates.
This motivated the development of several alternative conforming schemes to reduce the computational
overhead, e.g., the Bell element
as an modification of the Argyris element with less degrees of freedom \cite{ciarlet_finite_2002}.
In contrast, the non-conforming adaptive Morley FEM is known to be optimal and comes with an implementation in only 30
lines of MATLAB \cite{carstensen_discrete_2014}.

Adaptive mesh-refinement techniques for the many $C^1$ conforming FEMs remained unclear from the theoretical perspective until the preceding work
of Carstensen and Hu \cite{carstensen_hierarchical_2021}.
Their slight modification to the Argyris FEM comes with an adaptive algorithm and an efficient multilevel solver
at the cost of a negligibly increased computational effort. 
They prove optimal convergence rates 
and optimality of the proposed multilevel solver based on a multigrid V-cycle for the so-called hierarchical Argyris FEM.
Naturally, the comparison with the standard Argyris FEM is an important aspect from the practical viewpoint. 
This and numerical evidence of the optimality in physically relevant benchmarks justify a
rehabilitation of the Argyris element for fourth-order problems.
Since \cite{carstensen_hierarchical_2021} exclusively discusses homogeneous boundary conditions, the application to
meaningful models in solid mechanics requires the extension to general boundary conditions.
\subsection{Plate problem and FEM model}%
\label{sub:Model problem}
This paper considers the biharmonic plate equation with inhomogeneous mixed boundary conditions as a model example of a
fourth-order problem given by
\begin{align}
	\begin{aligned}\label{eqn:SF}
		\Delta^2 u= F && \text{in }\Omega,&&
		u=g&&\text{on }\GC\cup\GS, &&
		\Dn u=\Dn g&&\text{on }\GC.
	\end{aligned}
\end{align}
The function $u\in H^2(\Omega)$ describes the displacement of a thin structure or plate with mid-section $\Omega$ under the influence of
a force $F$.
Different boundary conditions model how the plate is hold in place, see the survey \cite{sweers_survey_2009}.
Clamped boundary conditions apply on $\GC\subset\overline\Omega$ and prescribe the displacement and bending of the plate in terms of the
globally defined boundary data $g\in H^2(\Omega)$.
On $\GS\subset\overline\Omega$, the plate is simply-supported and only its displacement is fixed.
This paper extends the a posteriori analysis of \cite{carstensen_hierarchical_2021} to general boundary conditions under
reasonable assumptions.
The first one is classical in the theory of plates \cite{brenner_mathematical_2008} and requires the
relatively open 
boundary components $\GC, \GS\subset\widebar\Omega$ of co-dimension one to ensure that the test space
$$V \coloneqq \{v\in H^2(\Omega)\ |\ v=\Dn v=0 \text{ on }\GC, v=0\text{ on }\GS\}$$
solely contains the trivial affine function, i.e., $V\cap P_1(\Omega)=\{0\}$.
Thus, the weak Hessian $D^2$ defines the bilinear energy form
$a(v,w) \coloneqq (D^2 v, D^2 w)_{L^2(\Omega)}$ for $v,w\in H^2(\Omega)$ that 
is positive definite on $V$ and induces the energy norm $\trb{\bullet}\coloneqq a(\bullet,
\bullet)^{1/2}$.
The weak form of the plate problem \eqref{eqn:SF} for given force $F\in V'$ and boundary data $g\in H^2(\Omega)$ seeks
the displacement $u\in g+ V$ defined by
\begin{align}\label{eqn:WP}
	a(u,v) &= F(v)&&\text{ for all }v\in V.
\end{align}
The standard (resp.\ hierarchical) Argyris FEM on a triangulation $\T$ seeks an approximation $\uh$ of
$u$ in the standard (resp.\ extended) Argyris space $\Argyris(\T)\subset P_5(\T) \cap H^2(\Omega)$ of
conforming piecewise quintic polynomials.
This requires a discrete analogon of the boundary data $g$ and a natural choice comes from the nodal interpolation
operator $\I:C^2(\overline\Omega)\to\Argyris(\T)$ if $g\in C^2(\overline\Omega)$ is sufficiently smooth.
The discrete approximation $\uh\in \I g + V(\T)$ to \eqref{eqn:WP} in the discrete test space $V(\T)\coloneqq V\cap
\Argyris(\T)$ solves 
\begin{align}\label{eqn:DWP}
	a(\uh,\vh) &= F(\vh)&&\text{ for all }\vh\in V(\T).
\end{align}
The a posteriori analysis in this paper assumes that $F\in V'$ is the sum of an $L^2$ contribution plus point forces.
The main result establishes optimal convergence rates of the error $\trb{u - \uh}$ and
oscillations in an adaptive algorithm for slightly more regular boundary data $g\in B$ in the space
\begin{align}\label{eqn:B}
	B\coloneqq\left\{v\in C^2(\overline{\Omega})\ \middle|\ \begin{array}{ll}
			\text{for all edges }\Gamma\text{ of }\partial\Omega,& v_{|\Gamma}\in H^3(\Gamma)\text{ and } (\Dn v)_{|\Gamma}\in
	H^2(\Gamma)\end{array}\!\!\!\right\}.
\end{align}
Note that this only imposes conditions on $g_{|\GCS}$ and $\Dn g_{|\GC}$ (e.g., replace $g$ by any element from $g+V$).
\subsection{Outline}%
\label{sub:Outline}
Section \ref{sec:The standard and hierarchical Argyris FEM} introduces some notation for the standard and hierarchical
Argyris FEM for general boundary conditions and the adaptive algorithm.
This paper is split into an analytical part preceding the numerical benchmarks in the second part.
The a posteriori error analysis in section \ref{sec:A posteriori analysis and optimality} extends
the optimal rates of the
hierarchical Argyris AFEM from \cite{carstensen_hierarchical_2021} to mixed inhomogeneous boundary conditions.
This directly leads to the equivalence of some computable a posteriori error estimator $\eta(\T)$ to the exact error (up to
oscillations) also for the standard Argyris FEM and 
motivates the comparison of the standard and hierarchical Argyris AEFM with the $\eta(\T)$ driven adaptive
algorithm in section \ref{sec:Numerical results}.
Section \ref{sec:Numerical2} discusses the application of multilevel-preconditioned iterative schemes for the solution
of the discrete problem \eqref{eqn:DWP}.
A reliable and efficient estimator of the algebraic error provides numerical evidence for the interoperability of multigrid
(MG) and preconditioned conjugated gradient (PCG) methods with the adaptive algorithm.
Section \ref{sub:Conclusion} concludes with some remarks.
\subsection{Overall notation}%
\label{sub:Overall notation}
Standard notation for Lebesgue and Sobolev spaces and their norms applies throughout this paper.
Let $H^s(K)$ abbreviate $H^s(\mathrm{int}(K))$ for closed $K\subset \R^2$.
Consider an open bounded Lipschitz domain $\Omega\subset\R^2$. 
The polygonal boundary $\partial\Omega$
with vertices $\VerticesO$ and edges $\EO$ decomposes into the relatively open, disjoint parts $\GC,\GS$ and into
$\partial\Omega\setminus(\GCS)$.
Let $P_k(K)$ denote the spaces of piecewise polynomials of total degree less than or equal $k\in\N_0$ on some triangle or edge $K\in\T\cup\E$
with diameter $h_K\in P_0(K)$.
The associated $L^2$ projection reads $\Pi_{K,k}:L^2(K)\to P_k(K)$ and is defined by the $L^2$ orthogonality
$(1-\Pi_{K, k})v\perp P_k(K)$ for all $v\in L^2(K)$.
Let
$$P_k(\T)\coloneqq\{p\in L^\infty(\Omega) : p_{|T}\in P_k(T)\text{ for all } T\in T\}$$
denote the space of piecewise polynomials on a triangulation $\T$ of the domain $\Omega$.
The partial derivatives $\partial_{v_1,\dots,v_j}^j$ for the $j\in\N_0$ directions $v_1,\dots,v_j\in\R^2$ define the
functional $\partial_{v_1,\dots,v_j}^j\delta_z:C^j(\overline\Omega)\to\R$ by
\begin{align*}
	\partial_{v_1,\dots,v_j}^j\delta_z(f) \coloneqq (-1)^j\delta_z(\partial_{v_1,\dots,v_j}^j f)
	=(-1)^j\partial_{v_1,\dots,v_j}^j f(z)&&\text{ for all }f\in C^j(\overline{\Omega})
\end{align*}
for the Dirac functional $\delta_z$ associated with some point $z\in\overline \Omega$.
\section{Adaptive standard and hierarchical Argyris FEM}%
\label{sec:The standard and hierarchical Argyris FEM}
This section defines the adaptive standard Argyris FEM, e.g., \cite{ciarlet_finite_2002}, and the adaptive hierarchical
Argyris FEM
\cite{carstensen_hierarchical_2021} for general boundary conditions.
\subsection{Triangulation}%
\label{sub:Triangulation}
Throughout this paper, $\T$ denotes a shape regular triangulation (in the sense of Ciarlet) of the polygonally bounded Lipschitz domain
$\Omega\subset\R^2$ with vertices $\mathcal{V}$ and edges $\E$ resolving the boundary conditions, i.e.,
$\bigcup \E(\Gamma_X) = \overline \Gamma_X$ for 
$\E(\Gamma_X)\coloneqq \{E\in\E\ :\ E\subset \overline\Gamma_X\}$ and $\Gamma_X=\GC, \GS$. 
The set $\E(\Omega)$ (resp.\ $\E(\partial\Omega)$) denotes the interior (resp.\ exterior) edges and the same notation
applies for the vertices $\Vertices$, edge-midpoints $\mathcal M\coloneqq\{\midE\ : \ E\in\E\}$, and nodes $\mathcal N\coloneqq
\Vertices\cup\mathcal{M}$.
Given a triangle $T\in\T$, denote the unit outer normal vector on the edges $E\in\E(T)$ of $T$ by $\nu_T$.
Associate every edge $E\in\E$ with a unit tangential $\tau_E$ and normal $\nu_E$ of fixed orientation.
If the context allows, the index $E$ with partial derivatives in directions $\tau_E, \nu_E$ is omitted.
The jump $[q]_E\in L^2(E)$ of $q\in H^1(\T)$ along an interior edge $E=T_+\cap T_-\in\E(\Omega)$ reads $[q]_E\coloneqq q_{|T_+}-
q_{|T_-}$ and $[q]_E\coloneqq q$ for boundary edges $E\in\E(\partial\Omega)$.

Fix an initial triangulation $\T_0$ of $\Omega$ with vertices $\Vertices_0$ and the set $\TT(\T_0)$ of all
admissible refinements generated by the newest-vertex bisection (NVB) \cite{stevenson_completion_2008, stein_finite_2017} of $\T_0$.
Note that the mesh-closure estimate requires no initial condition for $\T_0$ in two space dimensions
\cite[Thm.\ 2]{karkulik_2d_2013}.
For simplicity, some constructions related to a triangulation $\T\in\TT(\T_0)$ are formulated in terms of some sequence $(\T_0,
\T_1,\dots,\T_L=\T)$ of successive NVB refinements.
However, this construction will only depend on $\T_0$ and not on the chosen sequence.
\subsection{Standard and extended Argyris space}%
\label{sub:Discrete spaces}
The \textit{standard Argyris space} on $\T$, associated to the quintic Argyris element, consists locally of quintic polynomials and reads
\begin{align}\label{eqn:A_std}
	\Astd(\T)\coloneqq\left\{v_h\in P_5(\T)\cap C^1(\Omega) : D^2v_h\text{ is continuous at every } z\in\Vertices\right\}.
\end{align}
The extension to higher-order elements, e.g., the Argyris
element of order seven is straightforward and not addressed in this paper.
Notice that the continuity requirement of the Hessian at $z\in\Vertices$
makes the standard Argyris space not
hierarchical, i.e., in general $\Astd(\T)\not\subseteq \Astd(\Tt)$ for a refinement
$\Tt$ of $\T$.
In fact, the second-order normal-normal derivative $\partial_{\nu_E\nu_E}^2$ at
some edge's midpoint $z=\midE$ could be discontinuous across an edge $E\in\E(\Omega)$ in $\Astd(\T)$, whereas \eqref{eqn:A_std}
enforces its
continuity in $\Astd(\Tt)$ for any refinement $\Tt$ of $\T$ that contains $z\in\hat\Vertices$. %
The \textit{extended Argyris space} 
\begin{align}\label{eqn:A_ext}
	\Aext(\T) \coloneqq \Astd(\T_0) + \Astd(\T_1) + \dots + \Astd(\Tl)
\end{align}
is hierarchical and a minimal extension with respect to the sequence of successive refinements $\T_0, \dots, \Tl=\T\in\TT(\T_0)$.
The dependence on the initial triangulation is clear but the definition is in fact independent of the
sequence in \eqref{eqn:A_ext}, see
\cite{carstensen_hierarchical_2021} for further details.
Throughout this paper, let $\Argyris(\T)$ denote either $\Astd(\T)$ or $\Aext(\T)$ whenever no distinction is needed and
define the (conforming) discrete test space
\begin{align}\label{eqn:V}
	V(\T)\coloneqq\{v_h\in\Argyris(\T)\ :\ v_h=0\text{ on }\GCS\text{ and }\Dn v_h=0\text{ on }\GC\}=\Argyris(\T)\cap V.
\end{align}
\subsection{Local coordinate system}%
\label{sub:Local coordinate system}
The correct resolution of the boundary data and the degrees of freedom for the hierarchical Argyris FEM require control
of certain partial derivatives at the vertices.
Fix two directions $\{\tauz, \nuz\}$ for each vertex $z\in\Vertices$, spanning $\R^2$ (thought of as a local coordinate
system), under the two following conditions.
These are given in terms of some sequence $\T_0, \dots, \T_L=\T$ of successive refinements of the initial triangulation
$\T_0$ but only depend on $\TT(\T_0)$.
\begin{table}[]
	\centering
	\begin{tabular}{|c|c|c|c|c|c|}
		\hline
		\multicolumn{2}{|c|}{} & \multicolumn{2}{|c|}{$\omega=\pi$} &\multicolumn{2}{|c|}{$\omega\ne\pi$} \\\hline
		$E_0$ & $E_1$ & $\{\tauz, \nuz\}$ & $J(z)$ &$\{\tauz, \nuz\}$ & $J(z)$\\\hline
		$\GC$ & $\GC$ & $\{\tau_{0}, \nu_{0}\}$ & $\{1,2,3,4,5\}$ & $\{\tau_{0}, \nu_{0}\}$ & $\{1,2,3,4,5,6\}$\\
		$\GC$ & $\GS$ & $\{\tau_{0}, \nu_{0}\}$ & $\{1,2,3,4,5\}$ & $\{\tau_{0}, \nu_{0}\}$ & $\{1,2,3,4,5\}$\\
		$\GC$ & $\GF$ & $\{\tau_{0}, \nu_{0}\}$ & $\{1,2,3,4,5\}$ & $\{\tau_{0}, \nu_{0}\}$ & $\{1,2,3,4,5\}$\\
		$\GS$ & $\GS$ & $\{\tau_{0}, \nu_{0}\}$ & $\{1,2,4\}$ & $\{\tau_{0}, \tau_{1}\}$ & $\{1,2,3,4,6\}$\\
		$\GS$ & $\GF$ & $\{\tau_{0}, \nu_{0}\}$ & $\{1,2,4\}$ & $\{\tau_{0}, \nu_{0}\}$ & $\{1,2,4\}$\\
		$\GF$ & $\GF$ & any & $\emptyset$ & any & $\emptyset$\\\hline
	\end{tabular}
	\caption{Local coordinate system $\{\tauz, \nuz\}$ at
	$z\in\Vertices(\partial\Omega)$. Here, $\tau_{i}$ (resp.\ $\nu_i$) denotes the unit tangential (resp. normal)
	of $E_i\in\E(\partial\Omega), i=0,1$ with $z=E_0\cap E_1$ and $\omega$ the angle between $E_0, E_1$}
	\label{tab:tau_nu}
\end{table}
\newline
\textit{Condition 1:} If $z\in\Vertices(\Omega)\setminus\Vertices_0$ is a new interior vertex, then the NVB-algorithm yields $z=\midE$ for an edge $E$ of some previous
triangulation $\Tl, 0\leq \ell\leq L-1$. In this case, set $\tauz=\tau_E$ and $\nuz=\nu_E$.
\newline
\textit{Condition 2:} For a boundary vertex $z\in\Vertices(\partial\Omega)$, table
\ref{tab:tau_nu} provides a choice that depends on the boundary
conditions at the two boundary edges $E_0, E_1\in\E(\partial\Omega)$ that meet at $z=E_0\cap E_1$.\newline
No restrictions apply for the remaining cases where $z\in \mathcal{V}_0(\Omega)$ and the standard basis of $\R^2$ is a
natural choice.
\subsection{Degrees of freedom and nodal basis}%
\label{sub:Nodal basis and interpolation}
There is $m(z)\coloneqq 1$ degree of freedom (dof) $L_{z,1}\coloneqq\DnF\delta_z$ associated to each edge midpoint $z\in\mathcal M$.
This is the evaluation in the normal direction $\nu_E$ to the edge $E\in\E$ at $\midE=z$.
The other $m(z)\coloneqq 6$ (resp.\ $m(z)\in\{6,7\}$) dofs for the standard (resp.\ extended) Argyris space are associated with the
vertices $z\in\Vertices$ and consist of partial derivatives in the local coordinate system $\{\tauz, \nuz\}$.
For the standard Argyris space, they read 
\begin{align*}\delta_z,\; \partial_{\tauz}\delta_z,\;
	\partial_{\nuz}\delta_z,\; \partial_{\tauz\tauz}^2\delta_z,\; \partial_{\tauz \nuz}^2\delta_z,\; \partial_{\nuz
\nuz}^2\delta_z&&\text{ for }z\in\Vertices\end{align*}
and for the extended Argyris space
\begin{align*}
	&\delta_z,\; \partial_{\tauz}\delta_z,\;
	\partial_{\nuz}\delta_z,\; \partial_{\tauz\tauz}^2\delta_z,\; \partial_{\tauz \nuz}^2\delta_z,\; \partial_{\nuz
\nuz}^2\delta_z&&\text{ for }z\in\Vertices(\partial\Omega)\cup\Vertices_0,\\
			   &\delta_z,\; \partial_{\tauz}\delta_z,\;
			   \partial_{\nuz}\delta_z,\; \partial_{\tauz\tauz}^2\delta_z,\; \partial_{\tauz \nuz}^2\delta_z,\; \partial_{\nuz
			   \nuz}^2\delta_z^+,\; \partial_{\nuz
\nuz}^2\delta_z^-&&\text{ for }z\in\Vertices(\Omega)\setminus\Vertices_0.
\end{align*}
They are enumerated in this order as $L_{z,1}, ..., L_{z, m(z)}$.
Recall $\nuz=\nu_E$ for every new vertex $z\in\Vertices(\Omega)\setminus\Vertices_0$ where $z=\midE$ for some historical
edge $E$.
The modification for the extended Argyris space is a split of the normal-normal derivative evaluation $\partial_{\nuz
\nuz}^2\delta_z$($=\DnnF\delta_z$) at these vertices into the one-sided evaluations
\begin{align}\label{eqn:Dnn_Aext}
	\partial_{\nuz\nuz}^2\delta_z^{\pm}\equiv\DnnF\delta_z^\pm\coloneqq\lim_{x\in H_\pm(z)\to z}\DnnF\delta_x
\end{align}
in the half-planes 
$H_\pm(z)\coloneqq\{x\in\R^2: \pm(x-z)\cdot \nu_E\geq0\}$.
This allows $\DnnF v_h(z)$ to attain distinct values in $H_+(z)$ and $H_-(z)$ for $v_h\in\Aext(\T)$ at any such
vertex $z\in\mathcal V(\Omega)\setminus\mathcal V$. 
Indeed, this modification is enough to obtain hierarchical spaces and shows independence of the chosen sequence in the
definition of \eqref{eqn:A_ext}, see \cite{carstensen_hierarchical_2021} for further details.
Recall the set of nodes $\mathcal{N}\coloneqq \Vertices\cup\mathcal{M}$ and denote the unique nodal basis (dual to the
dofs)
 of $\Argyris(\T)$ by $\mathcal B\coloneqq\{\vfzj\ :\ z\in\Nodes, j=1,\dots,m(z)\}$.
 The choice of the local coordinates with $J(z)$ from table \ref{tab:tau_nu} for boundary vertices
 $z\in\Vertices(\partial\Omega)$ ensures that 
 \begin{align}\label{eqn:VT_basis}
	\{\vfzj\in\mathcal B\ :\ j\not\in J(z)\text{ for }z\in\Vertices(\partial\Omega)\text{ or }z\not\in
	\mathcal{M}(\GC)\}\subset\mathcal{B}
\end{align}
is a basis of the discrete test space $V(\T)$ from \eqref{eqn:V} as the following result shows.
\begin{prop}\label{prop:VT_characterization} With $J(z)$ for $z\in\Vertices(\partial\Omega)$ from table
	\ref{tab:tau_nu}, it holds that
\begin{align*}
	V(\T)=\left\{v_h\in\Argyris(\T)\ :\ \begin{aligned}
			L_{z,j}(v_h)&=0\text{ for all }
			z\in\Vertices(\partial\Omega),j\in J(z)\text{ and}\\ L_{z,1}(v_h)&=0\text{ for all }z\in\mathcal M(\GC)
		\end{aligned}\right\}.
\end{align*}
\end{prop}
\begin{proof}
Let $E=\mathrm{conv}\{P_0,
P_1\}\in\E(\partial\Omega)$ denote some boundary edge with normal $\nu\coloneqq\nu_E$ and tangential $\tau\coloneqq\tau_E$ and consider any $v_h\in\Argyris(\T)$.
It is well known, e.g., \cite{ciarlet_finite_2002}, that $v_{h|E}\equiv0$ vanishes if and only if the nodal value of
$v_h$ and its first two tangential derivatives along $E$ vanish at both endpoints $P_0$ and $P_1$, i.e., 
\begin{align}\label{eqn:v_h_0}\delta_z(v_h)=\Dt\delta_z(v_h)=\Dtt\delta_z(v_h)=0\quad\text{for } z=P_0, P_1.\end{align}
Similarly $(\Dn v_h)_{|E}\equiv0$ holds if and only if 
\begin{align}\label{eqn:Dn_v_h_0}\Dn\delta_z(v_h)=\Dtn\delta_z(v_h)=0\quad\text{for } z=P_0,
P_1\quad\text{and}\quad\Dn\delta_{\midF}(v_h)=0.\end{align}
With $J(z)$ and $\{\tauz, \nuz\}$ from table \ref{tab:tau_nu}, the conditions \eqref{eqn:v_h_0}--\eqref{eqn:Dn_v_h_0} translate into
equivalent assertions in terms of the dofs.
This shows the asserted identity.
\end{proof}
\noindent Notice the special treatment of a corner $z\in\VerticesO$ of the domain $\Omega$ between edges $E_0,E_1\in\E(\GS)$ in table
\ref{tab:tau_nu} where
the mixed derivative $\partial_{\tau_0 \tau_1}^2$ remains a degree of freedom in $V(\T)$.

\subsection{Interpolation of boundary data}%
\label{sub:Essential boundary conditions}
The duality relation between the dofs and the nodal basis defines
the nodal interpolation operator %
$\I:C^{2}(\overline{\Omega})\to\Argyris(\T),$
\begin{align}\label{eqn:I}
	\mathcal{I} v \coloneqq \sum^{}_{z\in\mathcal{N}} \sum^{m(z)}_{j= 1} \Lzj(v)\bzj&&\text{ for all }v\in
	C^{2}(\overline\Omega).
\end{align}
The following best-approximation property motivates the choice $\mathcal{I} g\in\Argyris(\T)$ for the discrete boundary
data $g\in B$ in the space of admissible boundary data $B$ from \eqref{eqn:B}.
\begin{lemma}[edge best-approximation]\label{lem:best_approx}
	Consider $v\in C^{2}(\overline{\Omega})$ and set $\cb\coloneqq(1-45\pi^{-4})^{-1/2}$.
	If $v_{|E}\in H^{3}(E)$ and $\Dn v_{|E}\in
	H^{2}(E)$ for some edge $E\in\E$, then
	\begin{align*}
		(a) &&&\left\Vert \Dttt (1-\I)v\right\Vert_{\LFace}^{} = \phantom{\cb}\left\Vert (1-\Pi_{E, 2})\Dttt
			v\right\Vert_{\LFace}^{},\\
		(b) &&&\left\Vert \Dttn (1-\I)v\right\Vert_{\LFace}^{} \leq \cb \left\Vert (1-\Pi_{E, 2})\Dttn
			v\right\Vert_{\LFace}^{}.
	\end{align*}
\end{lemma}
		\begin{proof}
			$(a)$ Repeated integration by parts and the exactness of the interpolation  $\I$ at the vertices
			show $\Dttt\I v = \Pi_{E, 2}\Dttt v$. Indeed, for arbitrary $p_2\in P_2(E)$, $\Dttt p_2\equiv 0$ verifies
			\begin{align*}
			\langle p_2, \Dttt(1-\I)v\rangle_{\LFace} &= \langle\Dttt p_2, (1-\I)v\rangle_{\LFace} = 0.
		\end{align*}
		$(b)$ Let $\bE\in P_2(E)$ denote the edge-bubble function on $E=\mathrm{conv}\{P_0,P_1\}$ that vanishes at both endpoints and attains
		$1=\bE(\midE)$ at the midpoint.
		Consider an arbitrary $p_2\in P_2(E)$ and set $\dE\coloneqq{\Pi_{E,0}(\Dnb(1-\I)v)}/(\Pi_{E,0}\bE^2)\in\R$. 
		Since $\Dn(1-\I)v$ and $\bE^2$ vanish to first order at the endpoints $P_0$ and $P_1$,
		\begin{align*}
			\langle p_2, \Dttn(1-\I)v - \dE\Dtt\bE^2\rangle_{\LFace} &= \langle\Dtt p_2, \Pi_{E,0}\left(\Dnb(1-\I)v -
			\dE\bE^2\right)\rangle_{\LFace}=0
		\end{align*}
		holds by the integration by parts formula and 
		$\dE\Dtt\bE^2 = \Pi_{E, 2}\Dttn v-\Dttn\I v\in P_2(E)$ follows. 
		A direct
		computation reveals $\Pi_{E,0}\bE^2 = 8/15$ and $\Dtttt \bE^2 =
		384|E|^{-4}$.
		This, integrating by parts twice, and the stability of the $L^2$ projection
		show
		\begin{equation}\notag%
			\begin{split}
			\left\Vert \Pi_{E, 2}\Dttn v- \Dttn\I v\right\Vert_{\LEdge}^{2}
				&=\langle \dE\Dtttt\bE^2, \Dnb (1-\I)v\rangle_{\LFace}%
			\leq 720|E|^{-4}\left\Vert \Dnb(1-\I)v\right\Vert_{\LFace}^{2}.%
			\end{split}
		\end{equation}
		By definition of the interpolation $\I$, $\Dnb(1-\I)v$ vanishes at both endpoints and at $\midE$.
		A split of the domain of integration $E$ into $E_j\coloneqq\mathrm{conv}\{P_j,\midE\},j=0,1$ allows the application of a Friedrichs inequality followed by
		a Poincar\'e inequality \cite{payne_optimal_1960} with known constant $|E_j|/\pi$ on $E_j,j=0,1$ separately.
This, $|E_0|=|E_1|=|E|/2$, and the previous estimate show $\Vert \Pi_{E, 2}\Dttn v- \Dttn\I
		v\Vert_{\LEdge}^{2}\leq 45\pi^{-4}\Vert \Dttn (1-\I)v\Vert_{\LFace}^{2}$.
	This and the
		Pythagoras Theorem prove
		\begin{align*}%
			\left\Vert \Dttn (1-\I)v\right\Vert_{\LFace}^{2}&\leq\left\Vert (1-\Pi_{E, 2})\Dttn v\right\Vert_{\LEdge}^{2} +
			45\pi^{-4}\left\Vert \Dttn (1-\I)v\right\Vert_{\LFace}^{2}.
		\end{align*}
		Since $45\pi^{-4}<1$, an absorption of the rightmost term on the left-hand side concludes the proof.
		\end{proof}
\noindent Consequently, the distance of the interpolation error $(1-\I)g$ of $g\in B$ to the test space $V$ is bounded by boundary oscillations
$\osc(\mathcal{S}, g)$ defined on a subset of edges $\mathcal{S}\subseteq\E$ by
\begin{align*}
	\osc^2(\mathcal S, g)&\coloneqq \sum^{}_{E\in\E(\GC)\cap\mathcal{S}}|E|^3\|(1-\Pi_{E,2})\Dttn g\|_{L^2(E)}^2
	+
	\sum^{}_{E\in\E(\GCS)\cap\mathcal{S}} |E|^3\|(1-\Pi_{E,2})\Dttt g\|_{L^2(E)}^2.
\end{align*}
		\begin{lemma}\label{lem:biharm_estimate}
			There exists a constant $\Cbiharm>0$ solely depending on $\Omega$ such that for any $v\in B$,
			$$\min_{w\in (1-\I)v+\V}\trb{ w}^2\leq \Cbiharm^2
			\osc^2(\F(\partial\Omega),v).$$
		\end{lemma}

		\begin{proof}
			It is straight-forward to verify that $t=(\varphi,\psi)\in \prod_{\Gamma\in\EO}(H^{3/2}(\Gamma)\times
			H^{1/2}(\Gamma))$, defined on each edge $\Gamma\in\EO$ of the polygonal domain $\Omega$ by
			\begin{align*}
				\varphi_{|\Gamma} = \begin{cases}{}
					(1-\I)v&\text{ on }(\GCS)\cap\Gamma,\\
					0&\text{ else}
				\end{cases}\quad\text{ and }\quad \psi_{|\Gamma}= \begin{cases}{}
					\Dnb((1-\I)v)&\text{ on }\GC\cap\Gamma,\\
					0&\text{ else,}
				\end{cases}
			\end{align*}
			belongs to the domain of the continuous right inverse of the trace map $\gamma_1:H^2(\Omega)\to
			\prod_{\Gamma\in\EO}(H^{3/2}(\Gamma)\times
			H^{1/2}(\Gamma))$ from \cite{grisvard_singularities_1992}.
			The extension $\hat t\in H^2(\Omega)$ of $t$ to the whole domain lies in $(1-\I)v+V$ and the boundedness of the
			right-inverse verifies
			$$\min_{w\in (1-\I)v+\V}\trb{ w}^2\leq \|\hat t\|_{H^2(\Omega)}^2\lesssim \sum_{\Gamma\in\EO}\big(
			\left\Vert \varphi\right\Vert_{H^{3/2}(\Gamma)}^{2} + \left\Vert
			\psi\right\Vert_{H^{1/2}(\Gamma)}^{2}\big).$$
			Since $\varphi$ (resp.\ $\psi$) vanishes up to second (resp.\ first) order at the vertices of an
			exterior edge $E\in\E(\partial\Omega)$, the 
			Gagliardo-Nierenberg inequality \cite[Thm.\ 1]{brezis_gagliardo-nirenberg_2018}
			in combination with repeated Friedrichs inequalities shows
			$$\min_{w\in (1-\I)v+\V}\trb{ w}^2\lesssim \sum_{E\in\E(\GCS)}|E|^{3}
				\left\Vert \Dttt\varphi\right\Vert_{L^2(E)}^{2} + \sum^{}_{E\in\E(\GC)} |E|^{3} \left\Vert
			\Dtt\psi\right\Vert_{L^2(E)}^{2}.$$
		The application of lemma \ref{lem:best_approx} concludes the proof.
		\end{proof}
\subsection{Adaptive algorithm}%
\label{sub:Adaptive algorithm and optimality}
Let the source $F\in V'$ be given by an $L^2$ contribution $f\in L^2(\Omega)$ and point forces. 
Assume that the initial mesh is compatible with the point forces in the sense that their support is a vertex of the initial
triangulation $z\in\Vertices_0$, i.e., there are $\beta_z\in\R$ such that
\begin{align}\label{eqn:F}
	F(v)&\coloneqq (f, v)_{L^2(\Omega)}+\sum_{z\in\Vertices_0} \beta_z v(z)&& \text{ for all }v\in V.
\end{align}
Given $T\in\T$ and the discrete solution $\uh\in\Argyris(\T)$ to \eqref{eqn:DWP}, the refinement indicator reads
	\begin{equation}
		\begin{aligned}
			\eta^2(\T, T) &=
			 |T|^2\|f - \BiL \uh\|_{\LO}^2+ \osc^2(\E(T), g)\\
						  &+ \sum^{}_{E\in\E(T)\setminus\E(\GC)}|T|^{1/2}\|\JF{\Dnn u_h}\|_{\LFace}^2%
						  +\sum_{E\in\E(T)\setminus\E(\GCS)}
						  |T|^{3/2}\|\JF{\Dttn u_h + \DLn {u_h}}\|_{\LFace}^2
		\end{aligned}\label{eqn:eta_def}
	\end{equation}
	and drives the adaptive AFEM algorithm \ref{alg:AFEM} for the standard and hierarchical Argyris FEM.
	\begin{algorithm}
	\caption{$\Argyris$-AFEM}
\label{alg:AFEM}
\begin{algorithmic}
	\State{\textbf{Input:}} Initial triangulation $\T_0$, bulk parameter $0<\theta<1$
	\State{\textbf{Solve}} the discrete problem \eqref{eqn:DWP} on $\Tl$ for $\ul\in \Argyris(\Tl)$
	\State \textbf{Compute} for all $T\in\T_\ell$ the local estimations $\eta(\Tl, T)$ from \eqref{eqn:eta_def}
    \State \textbf{Mark} minimal subset $\mathcal S_\ell\subset\T_\ell$ with
	\State $$\theta \sum^{}_{T\in\Tl} \eta^2(\Tl,T) \leq \sum^{}_{T\in\mathcal{S}_\ell}\eta^2(\Tl,T) $$
	\State \textbf{Refine} $\T_\ell\longrightarrow \T_{\ell+1}$ as smallest NVB refinement of $\T_\ell$ with $\mathcal{S}_\ell\subseteq\T_\ell\setminus\T_{\ell+1}$
	\State{\textbf{Output:}} Sequence of triangulations $\T_\ell$ and discrete solutions $u_\ell$
\end{algorithmic}
\end{algorithm}

\section{A posteriori analysis and optimality}%
\label{sec:A posteriori analysis and optimality}
This section proves the optimality of the adaptive hierarchical Argyris FEM ($\Aext$-AFEM) for possibly inhomogeneous boundary
data $g\in B$ and source terms of the form \eqref{eqn:F} including point forces.
It follows that the error estimator $\eta(\T)\coloneqq(\sum^{}_{T\in\T} \eta^2(\T,
T))^{1/2}$ is reliable and efficient up to the oscillations $\osc(\E(\partial\Omega), g)$ and $\osc(\T, f)\coloneqq \sum^{}_{T\in\T} \|h_T^2(1-\Pi_{T,0})f\|_{L^2(T)}$.
The axioms of adaptivity require the set $\TT(N)\coloneqq\{\T\in\TT(\T_0)\ :\ |\T| - |\T_0|\leq N\}$ for $N\in\N$ and
lead to the generalisation of \cite[Thm.~6]{carstensen_hierarchical_2021} following \cite{grasle_conforming_2022}.
\begin{theorem}[rate optimality of AFEM]\label{thm:rate_optimal}
	There exists $0<\Theta<1$ and for all $0<s<\infty$ some constant
	$\Lambda_{eq}>0$, only depending on $\TT(\T_0)$, on $\Theta$ and on $s$,
	such that the sequence of triangulations $(\Tl)_\ell$ and discrete solutions $(u_\ell)_\ell$ from the $\Aext$-AFEM algorithm with
	$\theta\leq\Theta$ satisfies
	\begin{align*}
		\sup_{\ell\in\N_0}
			&(1+|\T_\ell| - |\T_0|)^{s}\left(\trb{u-u_\ell} + \osc(\Tl, f) + \osc(\El(\partial\Omega), g)\right)\\
			&\leq \Lambda_{\rm
			eq}\sup_{N\in\N_0}(1+N)^{s}\min_{\T\in\TT(N)}\left(\trb{u-u_h} + \osc(\T, f) + \osc(\E(\partial\Omega), g)\right).
	\end{align*}
\end{theorem}
\noindent The proof employs the axioms of adaptivity \cite{stein_finite_2017,carstensen_axioms_2014,carstensen_axioms_2017}
and departs with a
quasi-interpolation operator defined on $H^2(\Omega)$ in the spirit of \cite{carstensen_hierarchical_2021} that interpolates exactly at the vertices.
Let $\omega(T)\coloneqq \mathrm{int}(\bigcup\{K\in\T\ :\ T\cap K\ne 0\})$ denote the layer-1 patch around $T\in\T$ and
$\gamma_1(v)\coloneqq(v_{|\partial\Omega},\Dn v_{|\partial\Omega})$ the trace map on $H^2(\Omega)$ from \cite{grisvard_singularities_1992}.
		\begin{theorem}[discrete quasi-interpolation]\label{thm:quasi-interpolant}
			There are constants $\capx, \cc>0$ exclusively depending on $\TT(\T_0)$ such that for any admissible refinement $\Tt$ of $\T\in\TT(\T_0)$ there
			exists a
			linear operator
			$J:H^2(\Omega)\to \Argyris(\T)$ satisfying, for any $v\in \HtwoO$ and $\vht\in \Argyris(\Tt)$,
			\begin{enumerate}[(a)]
				\item $J(V) = \V(\T)$ and $Jv(z) = v(z)$ for all $z\in\Vertices$.
				\item $\vht_{|T} = (J\vht)_{|T}$ for any $T\in\T\cap\Tt$.
				\item If $\gamma_1(v)\in \gamma_1(\Argyris(\T))$ then $\gamma_1(v) =\gamma_1(J v)$ (preservation of
					discrete boundary data).
				\item $| Jv|_{H^m(\Omega)}^{} \leq \cc|v|_{H^m(\Omega)}$ for $ m=0,1,2$
					(stability).
				\item $\sum^{2}_{m=0} h_T^{m-2}|(1-J)v|_{H^m(T)} \leq\capx |v|_{H^2(\omega(T))}$ for any $v\in V$
					(approximation property).
			\end{enumerate}
		\end{theorem}
		\begin{proof}
			Given any $v\in \HtwoO$ and the nodal basis $\mathcal{B}$ from subsection \ref{sub:Nodal basis and
		interpolation} and follow the lines of \cite[Thm.\ 2]{carstensen_hierarchical_2021}. Define 
			$$Jv = \sum^{}_{z\in\mathcal{N}}\sum^{m(z)}_{j=1}  M_{z,j}(v)\varphi_{z,j}\in\Argyris(\T)$$
			as the discrete function with coefficients $M_{z,j}(v)$.
			The choice $M_{z,1}\coloneqq L_{z,1}=\delta_z$ ensures exact interpolation at the vertices
			$z\in\Vertices$.
			The other functionals $M_{z,j}$ are chosen as in
			\cite{carstensen_hierarchical_2021} for interior nodes and as in \cite[Sec.\ 6]{girault_hermite_2002} for the extension to
			the boundary such that $(b)$--$(c)$ hold. See \cite[Sec.\ 7]{girault_hermite_2002} for the stability $(d)$ and the
			approximation property $(e)$; further details can also be found in \cite[Sec.~5]{grasle_conforming_2022} and are omitted
			here.
		\end{proof}
\newcommand{\Cbiharmd}{C_{\rm osc,d}}
\noindent This allows a discrete version of lemma \ref{lem:biharm_estimate}. Set $\Cbiharmd\coloneqq(1+\cb)\cc\Cbiharm$
for the constants $\cb, \Cbiharm, \cc$ from lemmas \ref{lem:best_approx}, \ref{lem:biharm_estimate}, and theorem
\ref{thm:quasi-interpolant}, respectively.
\begin{cor}\label{cor:biharm_estimate_discrete}
	Let $\Tt$ be an admissible refinement of $\T\in\TT(\T_0)$ and $\It:C^2(\overline\Omega)\to\Aext(\Tt)$ the associated
	nodal interpolation operator. Then
	$$\min_{\hat w \in (\It-\I)v+\VTt}\trb{ \hat w}^2\leq \Cbiharmd^2\;
	\osc^2(\F(\partial\Omega)\setminus\Ft(\partial\Omega), v).$$
\end{cor}
\begin{proof}
	The definition of the nodal interpolation \eqref{eqn:I} verifies $\I\It=\I$ and theorem \ref{thm:quasi-interpolant}
	provides 
	the quasi-interpolation $\hat J: H^2(\Omega)\to \Aext(\Tt)$ onto the fine space with
	$\hat J\big((1-\I)\It v + V\big) = (\It - \I) v + V(\Tt)$.
	This, theorem \ref{thm:quasi-interpolant} $(d)$, and lemma \ref{lem:biharm_estimate} provide
	$$\min_{\hat w \in (\It-\I)v+\VTt}\trb{ \hat w}^2\leq \cc^2\min_{w \in (1-\I)\It v+V}\trb{w}^2\leq\cc^2 \Cbiharm^2\;
	\osc^2(\F(\partial\Omega), \It v).$$
	Since $(\Dttt \It v)_{|E}, (\Dttn \It v)_{|E}\in P_2(E)$ is a quadratic polynomial on a unrefined edge $E\in \E(\partial
	\Omega)\cap \hat\E(\partial\Omega)$, the oscillation contribution on $E$ vanishes.
	On a refined edge $E\in\E(\partial\Omega)\setminus \hat\E(\partial\Omega)$, $\Vert (1-\Pi_{E, 2})\Dttn \It
	v\Vert_{\LEdge}^{}\leq \Vert \Dttn \It v-\Pi_{E, 2}\Dttn v\Vert_{\LEdge}^{}$, a triangle inequality, lemma
	\ref{lem:best_approx}, and the stability of $L^2$ projections result in
\begin{align*}
	\left\Vert (1-\Pi_{E, 2})\Dttn \It v\right\Vert_{\LEdge}^{} &\leq \left\Vert \Dttn
(1-\It)v\right\Vert_{\LEdge}^{}+\left\Vert (1-\Pi_{E, 2})\Dttn
v\right\Vert_{\LEdge}^{}\\
&\leq (1+\cb)\left\Vert (1-\Pi_{E, 2})\Dttn
v\right\Vert_{\LEdge}^{}.
\end{align*}
An analogous estimation for the other term in the oscillations shows $\osc^2(\E(\partial\Omega),\It
v)\leq(1+\cb)^2\osc^2(\E(\partial\Omega)\setminus\hat\E(\partial\Omega),v)$ and 
the claim follows.
\end{proof}

\subsection{Axioms of adaptivity}%
\label{sub:Axioms of adaptivity}
The axioms of adaptivity provide a framework for the proof of theorem \ref{thm:rate_optimal}. A key notion is the
distance 
$\delta(\T, \Tt)\coloneqq \trb{\uh - \uht}$
between two admissible triangulations $\T,\Tt\in\TT(\T_0)$ with discrete solutions $\uh$ and $\uht$ to \eqref{eqn:DWP}.
The remaining parts in this section discuss axioms \eqref{eqn:A1}--\eqref{eqn:A3} and \eqref{eqn:A4e} for the proof
of theorem
\ref{thm:rate_optimal} and require the nestedness of the extended Argyris space $\Aext$.
\begin{theorem}
			\label{thm:A1_A2}
			For any admissible refinement $\Tt$ of $\T\in\TT(\T_0)$, discrete stability and reduction hold with constants
			$\Lambda_1, \Lambda_2\in\R$ only depending on $\T_0$, i.e., 
			\begin{align}
				|\eta(\Tt,\Tt\cap\T) - \eta(\T,\Tt\cap\T)|\leq \Lambda_1\delta(\T,\Tt)\label{eqn:A1}\tag{A1},\\
				\eta(\Tt,\Tt\setminus\T)\leq 2^{-1/4}\eta(\T,\T\setminus\Tt) +
				\Lambda_2\delta(\T,\Tt)\tag{A2}\label{eqn:A2}.
			\end{align}
		\end{theorem}
\noindent The proof of this theorem uses standard arguments
\cite{stein_finite_2017,carstensen_axioms_2014,carstensen_axioms_2017} as for the case of homogeneous boundary conditions and is
therefore omitted.

		\begin{theorem}[discrete reliability]\label{thm:A3}
			A constant $\Lambda_3\in\R$ solely depending on $\T_0$ exists for the $\Aext$-AFEM such
			that for any admissible refinement $\hat\T$ of $\T\in\TT(\T_0)$,
			\begin{align}
				\label{eqn:A3}\tag{A3}
				\delta(\T,\hat\T)\leq\Lambda_3\;\eta(\T,\T\setminus\hat\T).
			\end{align}
		\end{theorem}
		\begin{proof}
			Let $\uh\in \Aext(\T)$ and $\uht \in\Aext(\Tt)$ solve \eqref{eqn:DWP} on the triangulations
			$\T$ and $\Tt$, respectively.
			The nestedness $V(\T)\subseteq V(\Tt)$ verifies that
			the error $\et\coloneqq\uht-\uh$ lies in $\et\in (\It-\I)g+\VTt$, where
			$\It:C^2(\overline\Omega)\to\Aext(\Tt)$ is the nodal interpolation onto the fine space.
			In general, $(\It-\I)g\not\in\VTt$ and the proof departs with the split of the error $\et=\wt+\wgt$ into a
			conforming part $\wt\in \VTt$ and
			\begin{align}\label{eqn:argmin}
				\wgt \coloneqq \operatorname*{argmin}_{\hat w\in(\It-\I)g+\VTt}\trb{\hat w}.
			\end{align}
			The characterisation \eqref{eqn:argmin} of $\wgt$ shows $a$-orthogonality to $\VTt$. Thus, the Galerkin
			property and $J:H^2(\Omega)\to\Aext(\T)$ from theorem \ref{thm:quasi-interpolant} with $J\wt\in \VT$
			provide
			\begin{align}\label{eqn:split}
				\delta(\Tt,\T)^2 = a(\et,\et) = L(\wt-J\wt) - a(\uh, \wt - J\wt) + a(\et, \wgt).
			\end{align}
			Repeated integration by parts with the abbreviation $\vt\coloneqq (1-J)\wt$ and \eqref{eqn:F} result in
			\begin{align*}
				L(\vt) - a(\uh, \vt)=&\;(f-\LL u_h,\vt)_{\LO} + \sum^{}_{z\in\Vertices}
				\beta_z\vt(z)\\
									 &+ \sum^{}_{T\in\T} \sum_{E\in\E(T)} \left( \left\langle \Dttn u_h+\Dn\Delta\uh, \vt\right\rangle_{\LFace} -
										 \left\langle\Dnn u_h,
									 \Dn{\vt}\right\rangle_{\LFace}\right).
			\end{align*}
			Recall the exactness of the quasi-interpolation at the vertices $z\in\Vertices$, theorem
			\ref{thm:quasi-interpolant} (a), to see that $\vt(z)=0$.
			The steps in the proof of \cite[Thm.\ 4]{carstensen_hierarchical_2021} for this setting
			consist of Cauchy and trace inequalities as well as the approximation properties of $J$ in theorem
			\ref{thm:quasi-interpolant} and
			show
				$|L(\vt) - a(\uh, \vt)|\lesssim
				\eta(\T,\T\setminus\Tt)\trb{\wt}$.
				No contributions arise from $T\in\T\cap\Tt$ due to $\vt=(1-J)\wt=0$ from theorem
				\ref{thm:quasi-interpolant} $(b)$.
			The Pythagoras Theorem 
			$\trb{\et}^2 = \trb{\wt}^2 + \trb{\wgt}^2$
			verifies $\trb{\wt}\leq\trb{\et}$.
			A Cauchy inequality and corollary \ref{cor:biharm_estimate_discrete} for the remaining term in
			\eqref{eqn:split} provide
			$$a(\et, \wgt)\leq \trb{\wgt}\trb{\et} \leq \Cbiharmd\;
			\osc(\F(\partial\Omega)\setminus\Ft(\partial\Omega), g)\trb{\et}.$$
			The combination of the given arguments with $\osc(\F(\partial\Omega)\setminus\Ft(\partial\Omega), g)\leq\eta(\T,\T\setminus\Tt)$ leads to
			$\delta(\T,\Tt)^2\lesssim \eta(\T,\T\setminus\Tt) \delta(\T,\Tt)$.
			This proves the existence of $\Lambda_3$.
		\end{proof}
		\noindent The axiom of quasi-orthogonality is an immediate consequence
		\cite{carstensen_axioms_2014,carstensen_axioms_2017} of its weakened version with an
		epsilon, axiom \eqref{eqn:A4e} together with axioms \eqref{eqn:A1}--\eqref{eqn:A2}.
		\begin{lemma}[quasi-orthogonality with $\varepsilon>0$]\label{lem:A4e}
			For all $\varepsilon>0$ there exists a constant $\Lambda_4(\varepsilon)\in\R$
			such that for every $m,n\in\N_0$ and the output sequence $(\Tl)_\ell$ of the $\Aext$-AFEM, 
			\begin{align}\label{eqn:A4e}\tag{A$4_\varepsilon$}
				\sum^{m+n}_{\ell=m} \delta^2(\Tlnext, \Tl) \leq\Lambda_4(\varepsilon)\;\eta^2(\T_m) + \varepsilon
			\sum^{m+n}_{\ell=m} \eta^2(\Tl).\end{align}
		\end{lemma}
			\begin{proof}%
				\newcommand{\ulnext}{u_{\ell+1}}
				\newcommand{\level}{\mathrm{level}}
				Denote the discrete solution on level $\ell$ by $\ul$ and the edges by $\El$.
				Furthermore, abbreviate the test spaces on each level by $\V_{\ell}\coloneqq\V(\Tl)$.
				The Galerkin orthogonality, a Cauchy inequality, and corollary
				\ref{cor:biharm_estimate_discrete} show for any $L\geq\ell$,
				\begin{align}\notag
					a(u_{L+1}-\ulnext, \ulnext - \ul) &= \min_{w_{\ell+1}\in \ul-\ulnext+V_{\ell+1}}
					a(u_{L+1}-\ulnext, w_{\ell+1})\\
							&\leq\Cbiharmd\;\delta(\T_{L+1},\T_{\ell+1})\; \osc(\F_{\ell}(\partial\Omega)\setminus\F_{\ell+1}(\partial\Omega),
						g).\label{eqn:A4e_eqn}
				\end{align}
				Note that each fine edge $\hat E\in \E_{\ell}(\partial\Omega), \ell\geq m+1$, is generated by
				$k$-times bisection of some coarse edge $E\in\E_m(\partial\Omega), k\in\N_0$.
				Hence, $h_{\hat E}=h_E/2^k$ and 
				\begin{align*}%
					h_{\hat E}^{3}\|(1-\Pi_{\hat E,2})\partial_{\t\t\bullet}^3g\|_{L^2(\hat E)}^2\leq2^{-3k}h_E^3
				\|(1-\Pi_{E,2})\partial_{\t\t\bullet}^3g\|_{L^2(\hat E)}^2\end{align*}%
				follows with the abbreviation $\partial_{\t\t\bullet}^3=\Dttt$ or $\Dttn$.
				Since each edge $E$ appears at most once in the $n$-fold collection of
				$\E_{\ell+1}(\partial\Omega)\setminus\E_\ell(\partial\Omega)$, $\ell=m,\dots,m+n-1$, this %
				leads to the geometric series
				$$\sum^{m+n-1}_{\ell=m} \osc^2(\F_\ell(\partial\Omega)\setminus\F_{\ell+1}(\partial\Omega), g)\leq \sum^{\infty}_{\ell=m}
				8^{(m-\ell)}\osc^2(\F_m(\partial\Omega), g) = \frac{8}{7}\osc^2(\F_m(\partial\Omega), g).$$
				Thus, the oscillations decay sufficiently fast so that \eqref{eqn:A3}, \eqref{eqn:A4e_eqn}, the generalised young inequality
				with $\alpha^2=\varepsilon\Cbiharmd^{-1}\Lambda_3^{-2}>0$, and $\osc^2(\F_m(\partial\Omega), g)\leq\eta^2(\T_m)$
				show
				\begin{align*}
					\sum^{m+n}_{\ell=m}& \delta^2(\Tlnext, \Tl) 
					= \delta^2(\T_{m+n+1}, \T_m) -2 \sum^{m+n-1}_{\ell=m}
					a(u_{n+m+1}-\ulnext, \ulnext - \ul)\\
						   &\leq\Lambda_3^2 \eta^2(\T_m) + \Cbiharmd \sum^{m+n-1}_{\ell =
							   m}\left(\alpha^2\delta^2(\T_{n+m+1},\T_{\ell+1}) +
						   \alpha^{-2}\osc^2(\F_{\ell}(\partial\Omega)\setminus\F_{\ell+1}(\partial\Omega), g)\right)\\
						   &\leq(1+8\Cbiharmd^2/(7\varepsilon))\Lambda_3^2\eta^2(\T_m) +
						   \varepsilon \sum^{m+n}_{\ell =
						   m+1}\eta^2(\Tl).
				\end{align*}
				This holds for all $\varepsilon>0$ and the claim follows with $\Lambda_4(\varepsilon)\coloneqq
				(1+8\Cbiharmd^2/(7\varepsilon))\Lambda_3^2$.
			\end{proof}
		\begin{proof}[Proof (of theorem \ref{thm:rate_optimal})]
			By \cite{carstensen_axioms_2014,
			carstensen_axioms_2017}, the axioms of adaptivity \eqref{eqn:A1}--\eqref{eqn:A3} and \eqref{eqn:A4e} first imply the axiom of
			quasi-orthogonality and then, for suficiently small $\Theta_0>0$,  optimality with
			respect to $\eta$, i.e.,
			\begin{align*}
		\sup_{\ell\in\N_0}
			(1+|\T_\ell| - |\T_0|)^{s}\;\eta(\Tl)
			\lesssim\sup_{N\in\N_0}(1+N)^{s}\min_{\T\in\TT(N)}\;\eta(\T).
			\end{align*}
			It remains to show equivalence of the estimator to the error up to oscillations
			\begin{align*}
				\eta(\T) + \osc(\T, f)\approx\trb{u-u_h} + \osc(\T, f) + \osc(\E(\partial\Omega), g).
			\end{align*}
			This follows (also for the standard Argyris FEM) as in the homogeneous case;
			the reliability is a small modification of theorem \ref{thm:A3} and only requires $V(\T)\subset V$.
			The efficiency estimate follows with standard bubble-function techniques \cite{verfurth_posteriori_2013}
			where no contributions from the point forces in $F$ from \eqref{eqn:F} arise, because these bubble-functions
			vanish at the vertices $z\in\Vertices$.
			Further details are omitted.
\end{proof}
\section{Numerical evidence for optimal convergence}%
\label{sec:Numerical results}
\FloatBarrier
This section uses the AFEM algorithm \ref{alg:AFEM} with an exact solver for a numerical comparison of the standard Argyris
AFEM ($\Astd$-AFEM) with the hierarchical Argyris AFEM ($\Aext$-AFEM).
Four benchmarks, see table \ref{tab:Problems_A}, employ varying singularities and boundary conditions on the Square
($\Omega=(0, 1)^2$), L-shape ($\Omega=(-1, 1)^2\setminus [0,1)^2$) and Slit ($\Omega=(0,1)^2\setminus(\{0\}\times[0,1))$).
\begin{table}[]
	\centering
	\begin{tabular}{c|c|c|c|c}
		Problem & domains & boundary data & BC & remarks\\\hline
		B1, B2	&	Square, L-shape	&	$g\equiv0$ & clamped	&	uniform load ($f\equiv1$)\\
		B3	&	Slit	&	$g\ne 0$ & clamped	&	exact solution known\\
		B4	&	L-shape	&	$g\ne 0$ & mixed	&	point load ($F=\delta_z$)
	\end{tabular}
	\caption{Overview of the problems}
	\label{tab:Problems_A}
\end{table}
\begin{figure}[]
	\centering
	{\scalebox{.95}{\includegraphics{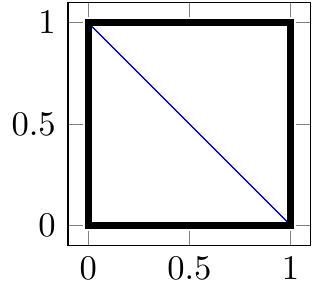}}
		\scalebox{.95}{\includegraphics{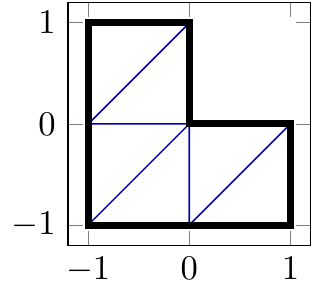}}
	\scalebox{.95}{\includegraphics{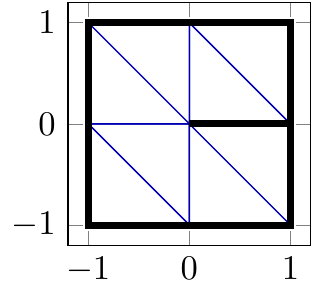}}
	\scalebox{.95}{\includegraphics{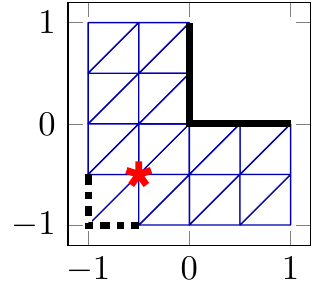}}}
	\caption{Initial triangulations $\T_0$ and boundary conditions (BC) $\GC$
		(\protect\includegraphics{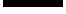}) and $\GS$ (\protect\includegraphics{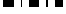}) from left to right: the Square, the L-shape, and the Slit with 
homogeneous BC and the L-shape with mixed BC and point force (\protect\includegraphics{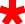})}
	\label{fig:E123_Domains}
\end{figure}

\subsection{Numerical realisation}%
\label{sub:Implementation}
\newcommand{\ldof}{\lambda}
Enumerate the nodal basis $\{\varphi_1, \dots,
\varphi_N\}$ of $V(\T)$ from \eqref{eqn:VT_basis} by $1$ to the number of degrees of freedom $N\coloneqq\mathrm{dim}(V(\T))$.
The algebraic formulation of \eqref{eqn:DWP} writes $\uhz=\sum_{j=1}^Nx_j\varphi_j\in V(\T)$ with coefficient vector $x\in
\R^N$ and seeks $\uhz=\uh-\I g$ such that
\begin{align}\label{eqn:ADWP}
	Ax = b
\end{align}
holds with the stiffness matrix
and right-hand side vector
\begin{align*}
	A &\coloneqq (a(\varphi_k,\varphi_j))_{k,j=1,\dots,N}\in \R^{N\times N}, &
	b &\coloneqq (F(\varphi_j) - a(\I g, \varphi_j))_{j=1,\dots,N}\in\R^N.
\end{align*}
The computation of $A$ and $b$ in the FEM fashion requires the evaluation of the $21$ local basis functions
 of the quintic Argyris finite element on each triangle $T\in\T$ at quadrature points.
 Possible approximation errors in the integration of non-polynomial expressions (e.g., from  $f$ from \eqref{eqn:F} and $g$) by quadrature are expected to be small and are
therefore neglected in this paper.
The transformation from the reference element 
\cite{kirby_general_2017,dominguez_algorithm_2008} allows an efficient evaluation of the local basis on the physical element.
Note that from the implementational viewpoint the hierarchical Argyris FEM only differs from the standard Argyris FEM in
that it treats the global $\partial_{\nuz, \nuz}^2\delta_z$ dof as the two degrees of freedom $\partial_{\nuz, \nuz}^2\delta_z^\pm$ from \eqref{eqn:Dnn_Aext}
for every $z\in\Vertices(\Omega)\setminus\Vertices_0$.
This section solves \eqref{eqn:ADWP} with the direct solver \texttt{mldivide} from the MATLAB standard library that is
behind the \texttt{\textbackslash} command.

\subsection{Benchmarks with homogeneous boundary conditions}%
\label{sub:Benchmarks with homogeneous boundary conditions}
This benchmark for the plate equation with uniform load $F\equiv 1$ consider homogeneous clamped boundary conditions
($g\equiv0$ and $\partial\Omega=\GC$) on
the Square ($\Omega=(0, 1)^2$) and the L-shape ($\Omega=(-1, 1)^2\setminus [0,1)^2$) with initial triangulations given in figure \ref{fig:E123_Domains}.
	Although the exact solution $u$ is unknown, the energy error of $u-\uh\in V$ can be computed by exploiting the Galerkin property for
	conforming discretisations, i.e.,
	$$\trb{u-u_h}^2 = \trb{u}^2-\trb{u_h}^2.$$
	The computation on a sufficiently fine mesh and multi-precision arithmetic led to the approximations
	$\trb{u}^2=3.8912007750677\times 10^{-4}$ %
	for the Square
	and $\trb{u}^2=3.57857007158618\times 10^{-3}$ %
	for the L-shape.
\begin{figure}[]
		\centering
		\includegraphics{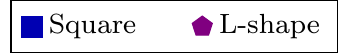}\\
		\includegraphics{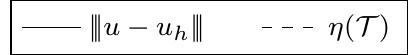}\\
		\hspace*{-.5em}\hbox{\includegraphics{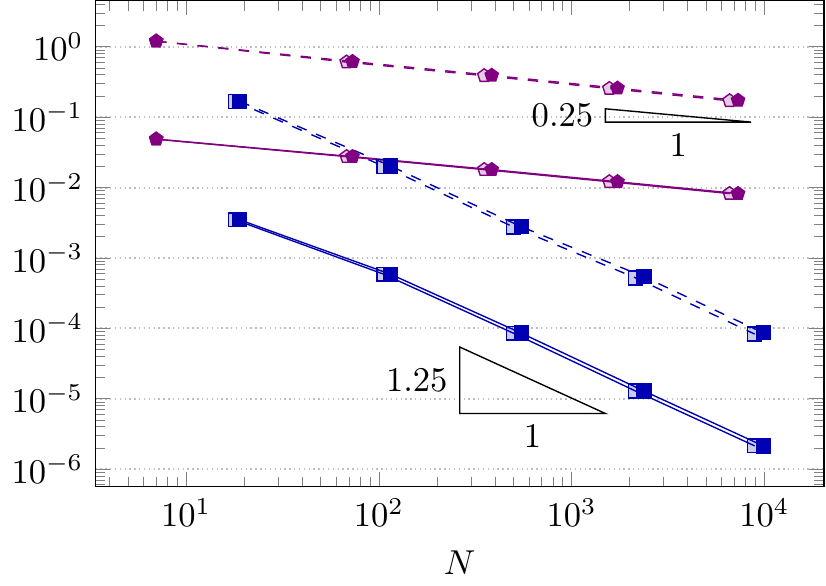}
		\includegraphics{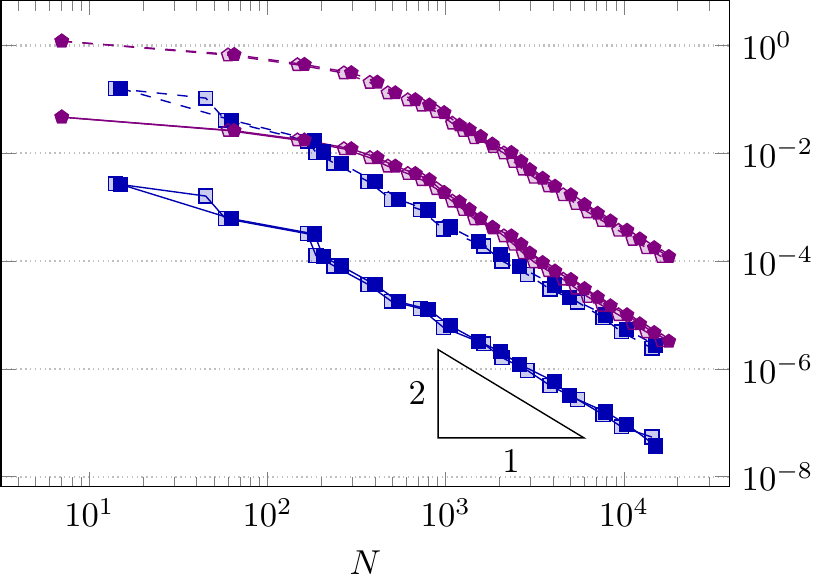}}
		\caption{Convergence history of the standard (opaque) and hierarchical (filled) Argyris AFEM on the L-shape and
		Square with uniform (left) and adaptive (right, $\theta=0.5$) mesh-refinement}
	\label{fig:E1_Convergence_uniform}
	\end{figure}
\begin{figure}[]
	\centering
	{\scalebox{1}{\includegraphics{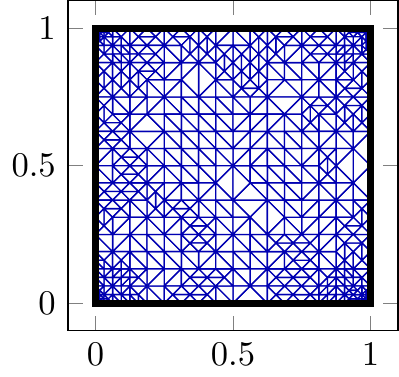}}
	\scalebox{1}{\includegraphics{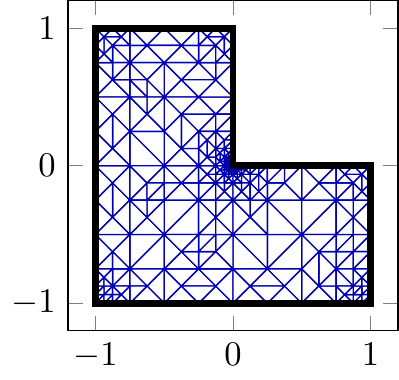}}}
	\caption{Adaptive triangulations $\T$ for $\theta=0.5$ of the $\Aext$-AFEM on the Square ($|\T|=861$) and 
		the L-shape ($|\T|=748$)}
	\label{fig:E1_Domains_adaptive}
\end{figure}

A uniform load on the L-shape with its corner singularity (at the origin) is the prototypical example of reduced
convergence rates for a uniformly refined
	mesh-sequence.
	Figure \ref{fig:E1_Convergence_uniform} does not only show an empirical suboptimal rate of $-1/4$ in the
	number of degrees of freedom \ndof\ for the L-shape but also a reduced rate of ${-5/4}$ on the
	Square. 
	This is shown for both the error in the energy norm $\trb{u-u_h}$ and the error estimator $\eta(\Tl)$ for the
	standard as well as the hierarchical Argyris FEM with uniform refinement.
	The reduced rates are due to corner singularities \cite{blum_boundary_1980} of the solution $u$ and underline the necessity of
	adaptive schemes even on the convex Square.

	Figure \ref{fig:E1_Domains_adaptive} shows refinement towards all (including convex) corners and a strong refinement
	towards the re-entering corner at the origin for the L-shape.
	Figure \ref{fig:E1_Convergence_uniform} also
	shows the theoretically predicted optimal convergence rates for the $\Aext$-AFEM  together with equivalence of the error
	estimator and the energy error.
	The same observations apply for the $\Astd$-AFEM.
	\FloatBarrier
	\subsection{Inhomogeneous boundary data on the Slit}%
	\label{sub:A benchmark with inhomogeneous BC}
	This benchmark problem considers the non-Lipschitz Slit ($\Omega=(0,1)^2\setminus(\{0\}\times[0,1))$) with pure
	clamped boundary $\GC=\partial\Omega$.
	The boundary data $g=u$ and source $F\equiv 1$ match the exact solution (in polar
	coordinates)
	\begin{align*}
		u(r,\varphi) &=-\frac{r^2}{16}\left(r^{1/2}\sin(\varphi/2) - \frac{r^2}{2}\sin^2(\varphi)\right).
	\end{align*}
\begin{figure}[ht]
		\centering
		\includegraphics{./Figures/Error_Legend_new_out.pdf}\\
		\hspace*{-3em}\includegraphics{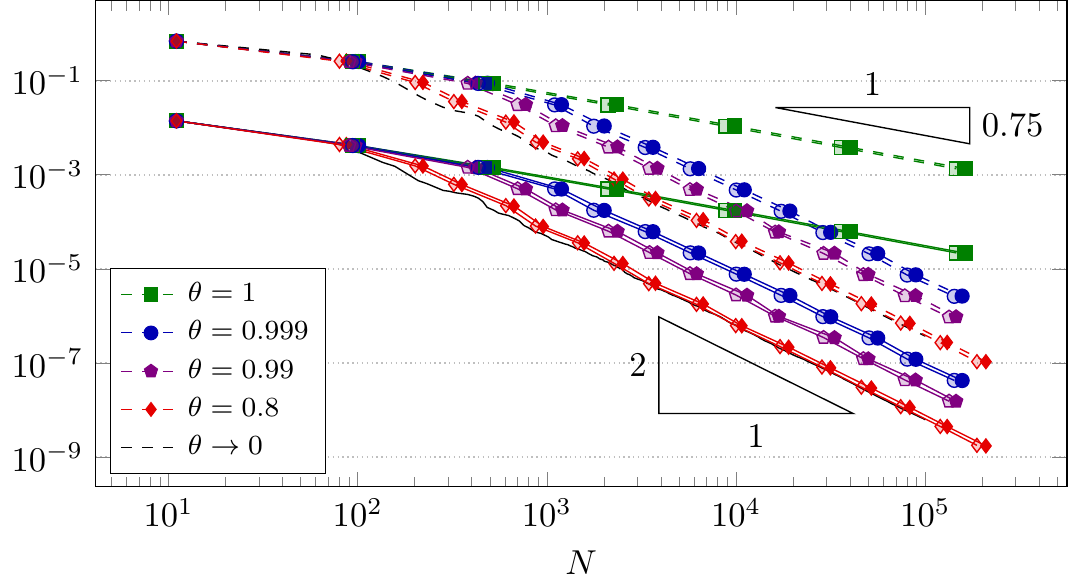}
			\caption{Convergence history of the standard (opaque) and hierarchical (filled) Argyris AFEM for different bulk
			parameter $\theta$}
		\label{fig:E2_Convergence}

	\end{figure}
\begin{figure}[]
	\centering
	{\scalebox{1}{\includegraphics{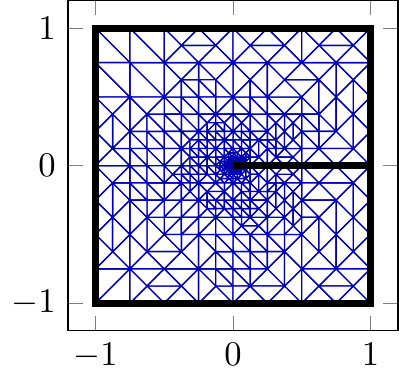}}
	\scalebox{1}{\includegraphics{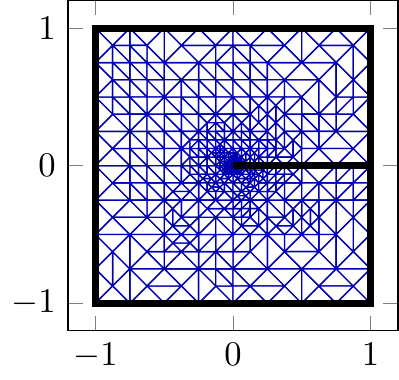}}}
	\caption{Adaptive triangulations $\T$ for $\theta=0.5$ of the $\Aext$-AFEM (left, $|\T|=700$) and of the
		$\Astd$-AFEM (right,
$|\T|=867$)}
	\label{fig:E2_Domain_adaptive}
\end{figure}
	\noindent Despite the singularity at the origin, the derivatives $\nabla u,D^2 u$ up to second order exist and vanish at the
	origin, hence the
	interpolation $\I
	u\in\Argyris(\T)$ is well defined.

	The theory on optimal convergence rates requires the bulk parameter $\theta$ to be sufficiently small.
	A first explicit computation \cite[Ex.\ 6.3]{carstensen_constants_2018} of the theoretical quantities for the Courant FEM 
 concludes optimality for $\theta\leq \Theta\coloneqq 2.6\times
	10^{-6}$.
	It is generally accepted that $\theta=0.5$ leads to optimal convergence in most practical scenarios.
	Figure \ref{fig:E2_Convergence} shows optimal rates even for $\theta$ close to one where $\theta=1$ abbreviates
	uniform refinement and $\theta\to0$ signals the opposite extreme with only $\mathrm{argmax}_{T\in\T}\eta(\T, T)$ marked
	for refinement.
	A higher bulk parameter $\theta$ initially leads to broader refinement and fewer steps of the adaptive algorithm in
	the pre-asymptotic regime.
	Yet, a choice of $0.5\leq \theta\leq 0.999$ produces similar iterations for $\ndof \geq
	10^3$.
	Even smaller values $\theta \leq 0.1$ lead to more refinement steps throughout but do not
	significantly improve on $\theta=0.8$.

	With the sole singularity of $u$ at the origin, figure \ref{fig:E2_Domain_adaptive} shows concentric refinement
	towards the origin.
	The adaptive mesh sequences obtained from $\Aext$-AFEM and $\Astd$-AFEM are qualitatively the same.
\subsection{Mixed boundary conditions and point load}%
\label{sub:A practical benchmark with multiple singularities}
In this benchmark the situation at hand is motivated by the biharmonic equation in the
context of plate bending. 
Consider a quadratic floor in some building, \eg skyscraper, made out of reinforced concrete. %
Suppose that the core of the building carries all the weight and that therefore the floor is embedded into the central
square.
Furthermore, there are supports around the outer corners that support the floor but do not fix tilting.
Since the layout is symmetric, the further considerations are reduced to the lower left quarter
so that the domain in consideration is represented by the L-shape with initial triangulation shown in figure~\ref{fig:E123_Domains}.
Displacement and bending are prescribed on $\GC=\{0\}\times[0,1)\cup [0,1)\times\{0\}$ whereas only displacement is fixed
on
$\GS=\{-1\}\times[-1/2,-1]\cup [-1,-1/2]\times\{-1\}$.

A point load $F=\delta_z\in H^{-2}(\Omega)$ at $z=(-1/2,-1/2)$ leads to a situation of no known solution $u$ to the
biharmonic equation
$$\Delta^2 u = \delta_z$$
with boundary data $g=10^{-3}\sin(\kappa x^3y^3\pi)$ introducing oscillations through $\kappa\geq0$.
	\begin{figure}[]
		\centering
		\includegraphics{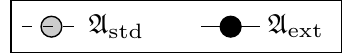}\\
		\hspace*{-3em}\includegraphics{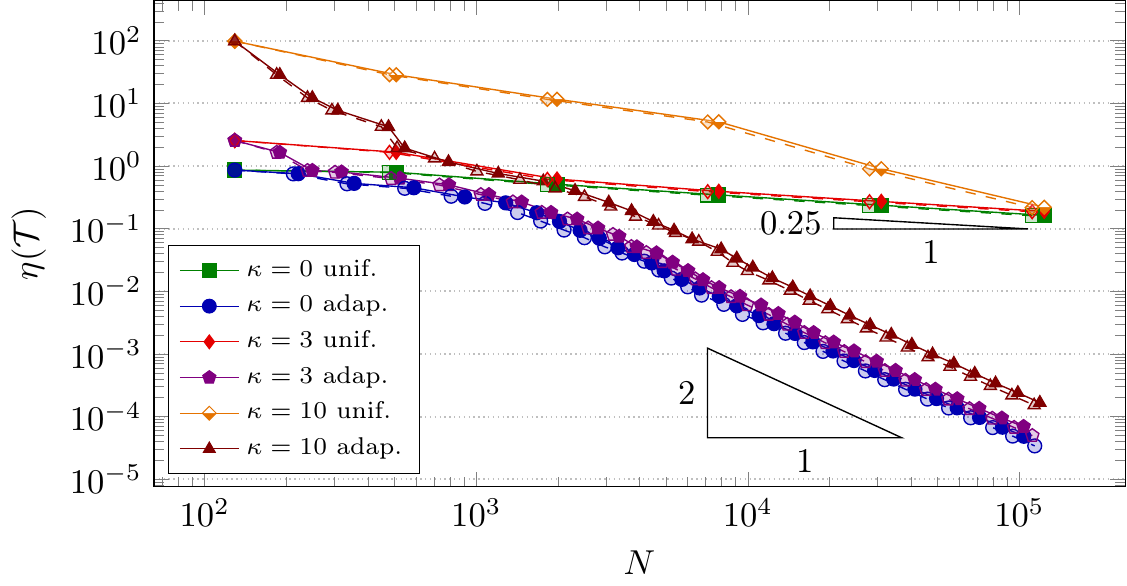}
		\caption{Convergence history for different values of $\kappa=0,3,10$}
		\label{fig:E4_Convergence}
	\end{figure}\noindent
\begin{figure}[]
		\centering
		\begin{minipage}{0.30\linewidth}
			\includegraphics{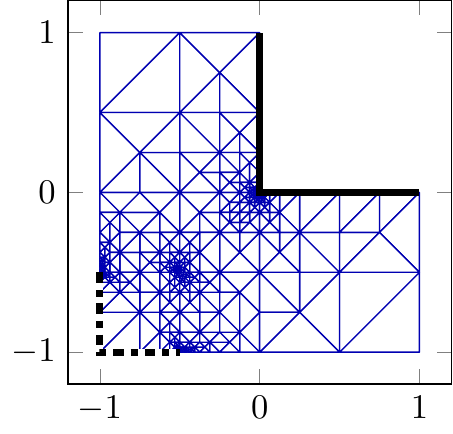}
		\end{minipage}
		\begin{minipage}{0.30\linewidth}
			\includegraphics{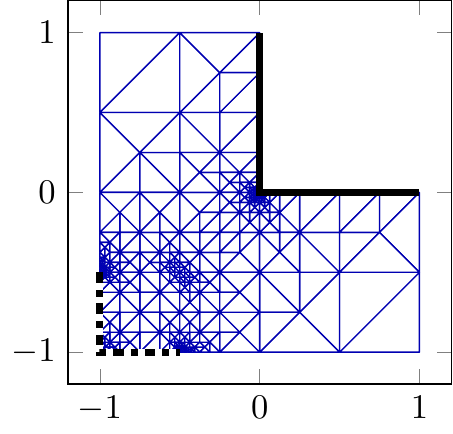}
		\end{minipage}
		\begin{minipage}{0.30\linewidth}
			\includegraphics{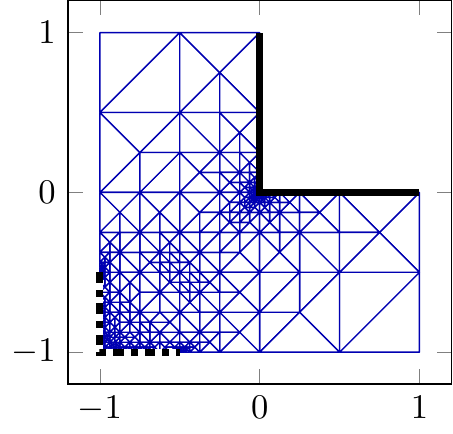}
		\end{minipage}
		\caption{Adaptive triangulations $\T$ from left to right; for $\kappa=0$ ($|\Vertices|=367$), for $\kappa=3$
			($|\Vertices|=372$), and for
		$\kappa=10$ ($|\Vertices|=340$)}
		\label{fig:E3_Meshes}
	\end{figure}
Mixed boundary conditions and the fact $F\not\in H^{-1}(\Omega)$ (but $F\in H^{-1-\varepsilon}(\Omega)$ for
$0<\varepsilon<1$) result in a problem with multiple singularities of different magnitude where the a priori mesh generation is not at all clear.
The adaptive algorithm leads to optimal rates for the standard and hierarchical Argyris AFEM, see figure
\ref{fig:E4_Convergence} for $\kappa=0,3,10$.
A large value of $\kappa$ introduces strong oscillations of the boundary data on $\GS$.
This initially leads to strong refinement on $\GS$ while the refinement for small $\kappa\geq0$ concentrates
on the location of the point force $z$, the origin, and the boundary of $\GS$.
Figure \ref{fig:E3_Meshes} shows these different intensities of the refinement for $\kappa=0,3,10$ on triangulations of
the $\Aext$-AFEM algorithm with $\theta=0.5$.
In fact, large parts, away from the singularities, were not further refined at all.
Anyhow, the oscillations of the boundary data are of higher order, see lemma \ref{lem:biharm_estimate}, so that their impact on uniform triangulations becomes
negligible for a small maximal mesh-size.

\section{AFEM with iterative multilevel solver}%
\label{sec:Numerical2}
This section compares the direct solver (from section \ref{sec:Numerical results}) with the iterative multigrid (MG) and preconditioned conjugated gradient (PCG)
methods for the solution to \eqref{eqn:DWP} of the hierarchical Argyris AFEM ($\Aext$-AFEM).
Throughout this section, the index $\ell$ refers to the level $\ell\in\N_0$ some
object is associated with.

\subsection{Adaptive multigrid V-cycle}%
\label{sub:Multigrid V-cycle}
Consider a sequence of successive refinements $(\Tl)_\ell$ with stiffness matrix $A_\ell$ and right-hand side vector and
$b_\ell$ at level $\ell\in\N_0$.
Recall the nodal basis $\{\vflj\}_{\ell=1}^{N_\ell}$ for the discrete test space $V(\Tl)$ of dimension
$N_\ell\coloneqq\mathrm{dim}(V(\Tl))$ and the algebraic formulation $A_\ell x_\ell = b_\ell$ of \eqref{eqn:DWP} from subsection \ref{sub:Implementation}.
Multigrid methods make use of the whole sequence of discretisations to solve 
$A_\ell x_\ell = b_\ell$.
The matrix version \cite{bramble_multigrid_1993} of the multigrid V-cycle for the hierarchical Argyris FEM \cite[Sec.\ 5]{carstensen_hierarchical_2021} requires the prolongation
matrix $P_\ell\in \R^{N_{\ell}\times N_\ell-1}$ that expresses a coarse function in terms of the fine basis functions,
i.e., $\varphi_k^{\ell-1} = \sum_{j=1}^{N_\ell} P_{\ell,jk}\vflj$ for $1\leq k\leq N_{\ell-1}$,
and the matrix $\opSl\in \R^{N_\ell\times N_\ell}$.
Let $I_\ell\coloneqq\{j\ :\ \vflj \not\in V(\T_{\ell-1})\}$ denote the indices of the new basis functions at level
$\ell$.
The \textit{local Gauß-Seidel} smoother $\opSl$ acts on the $I_{\ell}$-components of a vector $y_\ell\in \R^{N_\ell}$
by 
\begin{align}\label{eqn:Sl}(\opSl y_\ell)_j \coloneqq \begin{cases}{}
	(\tilde A_\ell ^{-1}y_{\ell|I_\ell})_j&\text{ if }j\in I_\ell,\\
	0&\text{ else,}
\end{cases}\end{align}
where $\tilde A_\ell=\mathrm{tril}(A_{\ell|I_\ell I_\ell})$ denotes the lower triangular part (including the diagonal)
of the submatrix
$A_{\ell|I_\ell I_\ell}=(A_{\ell, kj})_{k,j\in I_\ell}$ and $y_{\ell|I_\ell} = (y_{\ell, j})_{j\in I_\ell}$, see
\cite{bramble_analysis_1992} for the relation with the operator notation of $\opSl$ in
\cite{carstensen_hierarchical_2021}.
This way, $\opSl$ only acts on the components that correspond to new basis functions (either
associated to a new node $z\in\Nodes_\ell\setminus\Nodes_{\ell-1}$ or with support on a refined triangle
$T\in\T_\ell\setminus\T_{\ell-1}$).
\begin{algorithm}
\begin{algorithmic}
\State\textbf{Input:} $y_\ell\in V(\T_\ell), r\in\N$
	\If{$\ell=0$}
	\State \textit{Exact solve}: $\opB_0 y_0=\opA_0^{-1}y_0$
	\Else{}
	\State \textit{Pre-smoothing}: 
	$w_{j+1}\coloneqq w_j+ \opS_\ell(y_\ell-\opA_\ell w_j)$ for $w_0\coloneqq0$ and $j=0, \dots, r-1$
	\State \textit{Coarse-grid correction}: $w_{r+1} \coloneqq w_r+ P_\ell\opB\ellminus P_\ell^\top(y_\ell-\opA_\ell w_r)$
	\State \textit{Post-smoothing}: $w_{j+1}\coloneqq w_j+ \opS_\ell^\top(y_\ell-\opA_\ell w_j)$ for $j=r+1,\dots,2r$
	\State Set $\opB_\ell y_\ell\coloneqq w_{2r+1}$.
	\EndIf
 \State\textbf{Output:} $\opB_\ell y_\ell$
\end{algorithmic}
	\caption{(V$(r)$-cycle)}
\label{alg:v-cycle}
\end{algorithm}
The standard symmetric multigrid V-cycle \cite{bramble_multigrid_1993}, algorithm \ref{alg:v-cycle}, with $r$ pre- and
post-smoothing steps defines a uniform approximative inverse $\opBl$ of $\opAl$ 
in the spectral energy norm $\trb{M}_2\coloneqq\sup_{y_\ell\in \R^{N_\ell}\setminus 0}(y_\ell^\top A_\ell M
y_\ell)/(y_\ell^\top A_\ell y_\ell)$ for $M\in
\R^{N_\ell\times N_\ell}$.
The original proof for the $V(1)$-cycle holds for general $V(r)$-cycles.
\begin{theorem}[{\cite[Thm.\ 7]{carstensen_hierarchical_2021}}]\label{thm:uniform_approx}
For the $\Aext$-AFEM, there exists $c_\infty\in\R$ with
$$\sup_{\ell\in\N_0}\trb{\opI-\opBl\opAl}_2\leq \frac{c_\infty}{1+c_\infty}<1.$$
\end{theorem}
\begin{proof}
	Notice, e.g., from  \cite[Sec.\ 10]{bramble_multigrid_1993} and \cite[Sec.\ 3]{grasle_conforming_2022}, that $\opSl$ and $\opQl$ are the matrix representations of the local Gauß-Seidel relaxation operator in
	\cite{carstensen_hierarchical_2021} and the
	$L^2$ projection $V(\Tl)\to V(\T_{\ell-1})$. 
	This establishes $\opIBAl$ as the matrix representation of the operator $I-\mathcal{B}_\ell\mathcal{A}_\ell$ from
	\cite[Sec.\ 7.6]{carstensen_hierarchical_2021}. 
	The proof of \cite[Thm.~7]{carstensen_hierarchical_2021} provides $I-\mathcal{B}_\ell\mathcal{A}_\ell=\mathcal R^*\mathcal R$ with $\mathcal R=(I-\mathcal P_0)\prod_{j=1}^\ell
	(I-\mathcal
	Q_j)^r$ for $r$ smoothing steps and operators $\mathcal P_0, \mathcal Q_j, j=1,\dots,\ell$.
	The application of \cite{xu_method_2002} provides $c_\ell^{(r)}\in\R$ with \begin{align}
		\|I-\mathcal{B}_\ell\mathcal{A}_\ell\|_{\mathcal
	A_\ell}\coloneqq\sup_{\vl\in V(\Tl)\setminus
0}\frac{a((I-\mathcal{B}_\ell\mathcal{A}_\ell)\vl,\vl)}{a(\vl,\vl)}=\|\mathcal
	R\|_{\mathcal{A}_\ell}^2=\frac{c_\ell^{(r)}}{1+c_\ell^{(r)}}.
	\end{align} The explicit characterisation of
	$c_\ell^{(r)}\in\R$ shows $c_\ell^{(r)}\leq c_\ell^{(1)}\leq c_\infty$ with the uniform bound $c_\infty\in\R$
	\cite[Thm.\ 7]{carstensen_hierarchical_2021} and concludes the proof for general $r\in\N$.
\end{proof}

\subsection{Stopping criterion}%
\label{sub:Stopping criterion}
Let $x_\ell\in \R^{N_\ell}$ denote the exact solution to \eqref{eqn:ADWP}, i.e., $\ul = \ula + \I_\ell g$ solves
\eqref{eqn:DWP} for $\ula=\sum_{j=1}^{N_\ell} x_j\vflj$.
For any $\tilde u_\ell=\tula + \Il g$ with coefficient vector $\tilde x_\ell\in \R^{N_\ell}$ of $\tula\in V(\Tl)$, 
a reliable and efficient estimator 
$$\etaalg(\Tl, \tilde x_\ell)
\coloneqq\left((b_\ell - A_\ell\tilde x_\ell)^\top B_\ell (b_\ell - A_\ell\tilde x_\ell)\right)^{1/2}$$
for the algebraic error $\trb{\ul -\tilde u_\ell}$ comes from the approximate inverse $\opBl$
of the multigrid V-cycle.
This motivates the stopping criterion with tolerance $0<tol$ for initial solution $\txl^0\in\R^{N_\ell}$,
\begin{align}\label{eqn:stopping_c}
	\etaalg(\Tl, \tilde x_\ell) < tol\; \etaalg(\Tl, \tilde x_{\ell}^0).
\end{align}
\begin{lemma}[{\cite[Lem.~3.5]{grasle_conforming_2022}}]\label{lem:eta_alg_eq}
	If  $\trb{
	\opI-\opBl\opAl}_{2}^{}\leq C<1$, then for any $\tula\in V(\Tl)$ with coefficient vector
	$\txl\in\R^{N_\ell}$ and $\tilde u_\ell\coloneqq \tula+\I_\ell g\in\Aext(\T)$,
	\begin{align*}
		(1- C)^{1/2}\trb{ u_{\ell}-\tilde u_{\ell}}\leq \etaalg(\Tl, \tilde x_\ell)\leq (1+C)^{1/2}\trb{
			u_{\ell}-\tilde u_{\ell}}.
	\end{align*}
	
\end{lemma}
\begin{proof}
	With the error $e_\ell \coloneqq x_\ell-\tilde x_\ell$, the residual reads
	$A_\ell e_\ell = b_\ell - A_\ell \tilde x_\ell$ and 
	\begin{align}\label{eqn:error_B}\trb{u_{\ell}-\tilde u_{\ell}}^{2}\coloneqq a(u_{\ell,0}-\tilde u_{\ell,0}, u_{\ell,0}-\tilde u_{\ell,0}) =
	e_\ell^\top A_\ell e_\ell = e_\ell^\top A_\ell (I-B_\ell A_\ell)e_\ell + e_\ell^\top A_\ell B_\ell A_\ell e_\ell.\end{align}
	The definition of the norm, $\trb{
	\opI-\opBl\opAl}_{2}^{}\leq C$, and $\etaalg(\Tl,\tilde x_\ell)=e_\ell^\top A_\ell B_\ell A_\ell
	e_\ell$ show
	\begin{align*}
		\trb{u_{\ell}-\tilde u_{\ell}}^{2}
		&\leq C\;\trb{u_{\ell}-\tilde u_{\ell}}^{2}+\etaalg(\Tl,\tilde x_\ell).
	\end{align*}
	This proves the first inequality and the second follows by similar arguments after rearranging \eqref{eqn:error_B}.
\end{proof}
\noindent This section applies the (full) multigrid method (MG) with iterations
$$\tilde x_\ell^{j+1} = \tilde x_\ell^j + \opBl(b_\ell-\opAl \tilde x_\ell^j)$$
and the preconditioned conjugated gradient method (PCG) with preconditioner $\opBl$ as an example for Krylov subspace
methods \cite{saad_iterative_2003}.
The stopping criterion reads \eqref{eqn:stopping_c} and the initial solution $\tilde x_\ell^0$ is the coefficient vector of the solution at the previous level $
u_{\ell-1,0}$ for $\ell\geq 1$ or $\tilde x_0^0\coloneqq 0$; this is also known as nested iterations.
\begin{rem}[alternative stopping criterion]
	Gantner et al.~\cite{gantner_rate_2021} discuss an alternative stopping criterion that requires the evaluation of the error
	estimator $\eta(\Tl)$ at $\tilde u_\ell$ for each iteration and prove optimal rates of the overall adaptive algorithm
	(given a small enough tolerance). 
	In practice however, the evaluation of the error estimator $\eta(\Tl)$ from \eqref{eqn:eta_def} by
	quadrature dominates the computational time in the solve step of algorithm \ref{alg:AFEM}.
\end{rem}
The computation of the stopping criterion \eqref{eqn:stopping_c} in this paper is remarkably simple as the MG-update $B_\ell (b_\ell - A_\ell\tilde x_\ell)$ and
the residual $b_\ell - A_\ell\tilde x_\ell$ are already computed quantities of each MG or PCG iteration.

\subsection{Adaptive algorithm with inexact solve}%
\begin{figure}[]
		\centering
		\includegraphics{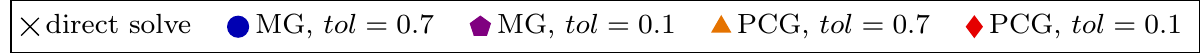}\\
		\hspace*{-1.5em}\hbox{\scalebox{1}{\includegraphics{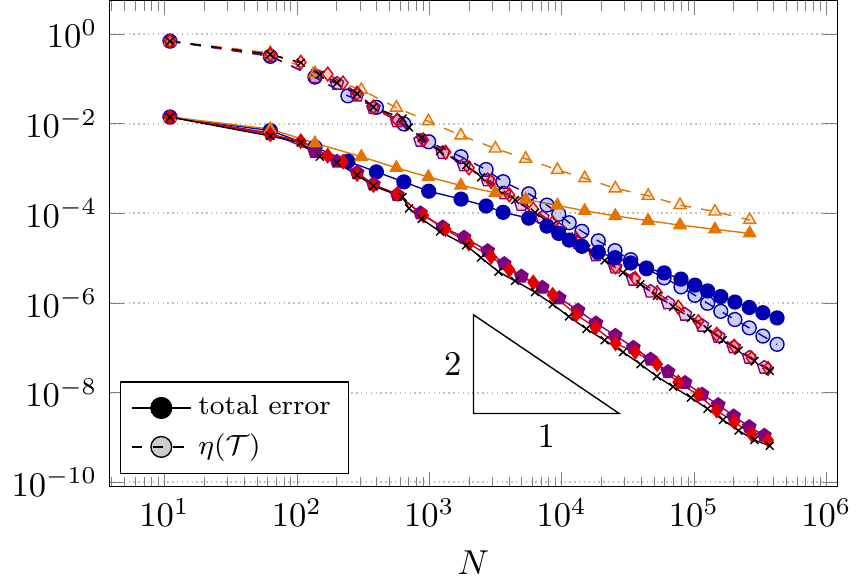}}
		\includegraphics{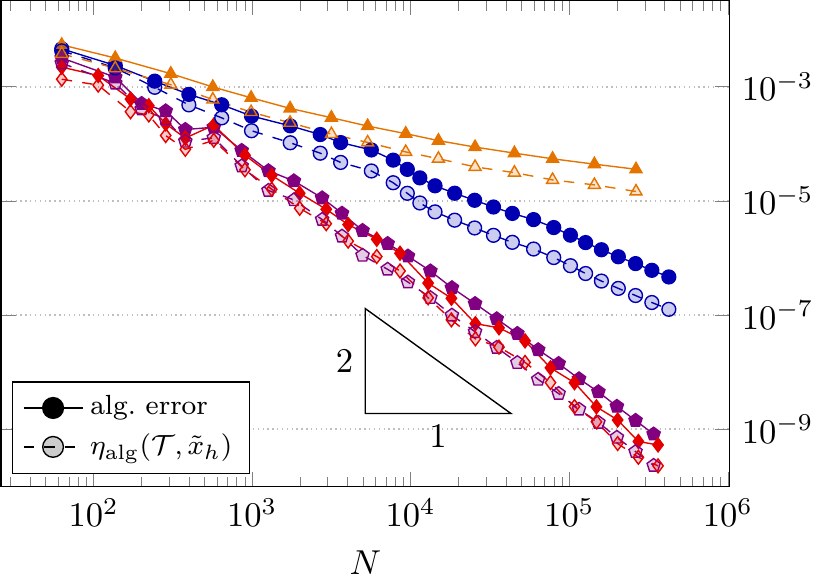}}

		\caption{Convergence history of the total error $\trb{u-\tilde u_h}$ (solid) and $\eta(\T)$ (dashed) for AFEM ($\theta=0.8$) using
		iterative multigrid (MG) and
		preconditioned CG (PCG) solvers compared with the direct solver for the Slit benchmark (B3)}
		\label{fig:Iter_Tol}
	\end{figure}

\label{sub:Adaptive algorithm with inexact solve}
A crucial ingredient in the proof of reliability of the error estimator $\eta(\T)$ is the Galerkin property
\begin{align}\label{eqn:Galerkin}
	\trb{u-\vh}^2 = \trb{u-\uh}^2 + \trb{\uh-\vh}^2
\end{align}
of the exact
discrete
solution $\uh$ to \eqref{eqn:DWP} and any $\vh\in \I g+V(\T)$.
For an approximation $\tuh\in\I g+V(\T)$ to $\uh$ this generally does not hold but shows that the
discretisation error $u-\uh$ is $a$-orthogonal to the algebraic error $\uh - \tuh$.
Since $\eta(\T)$ depends continuously on $\uh$, a small algebraic error $\trb{\uh- \tuh}$
is acceptable and leads to optimal rates of the $\Aext$-AFEM algorithm \ref{alg:AFEM} with inexact solve, see figure
\ref{fig:Iter_Tol} for the Slit benchmark B3 from subsection \ref{sub:A benchmark with inhomogeneous BC}.
For too coarse approximations ($tol=0.7$ in fig.\ \ref{fig:Iter_Tol}), the algebraic error does not converge with optimal rate and, thus by
\eqref{eqn:Galerkin}, the same holds true for the total error $\trb{u-\tuh}$.
In this case, the error estimator $\eta(\T)$, evaluated for $\tuh$, is not reliable and suggests a
better convergence rate.
Conversely, further (undisplayed) experiments with the benchmarks from section \ref{sec:Numerical results}
suggest that an optimal convergence of the algebraic error $\trb{\uh-\tuh}$ is also sufficient for reliability of $\eta(\T)$.
Hence, $\etaalg(\T, \tilde x_h)$ (equivalent to $\trb{\uh-\tuh}$
by lemma \ref{lem:eta_alg_eq} and theorem \ref{thm:uniform_approx}) serves as an indicator for optimal rates. 

Undisplayed numerical experiments with the related adaptive scheme from \cite{gantner_rate_2021} show a very similar
convergence behaviour and dependence on the tolerance of the stopping criterion compared to \eqref{eqn:stopping_c}.
This hints at a possible equivalence of both algorithms and motivates further investigations regarding a convergence
analysis extending \cite{gantner_rate_2021} to inhomogeneous boundary data.
	\begin{figure}[]
		\centering
		\includegraphics{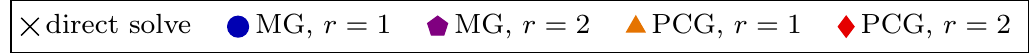}\\
		\hspace*{-2em}\hbox{\scalebox{1}{\includegraphics{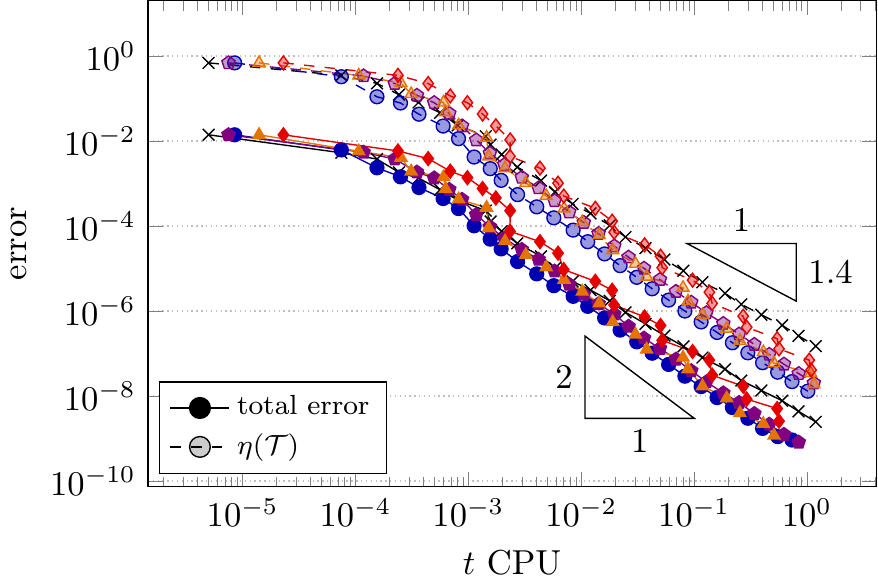}}
		\scalebox{1}{\includegraphics{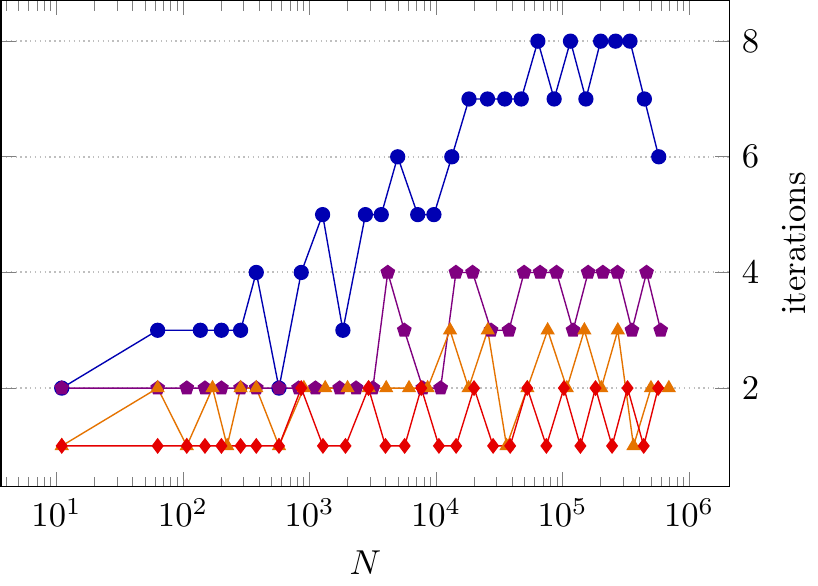}}}
		\caption{Convergence history of the total error $\trb{u-\tilde u_h}$ (solid) and $\eta(\T)$ (dashed) against the single-core CPU time $t$ of one
			MG/PCG iteration in seconds on Intel\textsuperscript{\textregistered}
			Xeon\textsuperscript{\textregistered} Gold 5222 CPU at 3.80GHz with 1000GB RAM (left) and the number of iterations
	(right) with $tol=0.1$ for the Slit benchmark (B3)}
		\label{fig:Timing}
	\end{figure}
	\subsection{Linear time complexity AFEM}%
	\label{sub:Optimal complexity AFEM}
	 Algorithm \ref{alg:AFEM} consists of the steps Solve, Estimate, Mark and Refine.
	 Each of these steps has a linear time complexity (in the number of degrees of freedom $N$)
	 if an appropriate solver
	 is applied in the Solve step, see \cite{pfeiler_dorfler_2020} for the linear complexity of Dörfler marking with
	 minimal cardinality.
	 Direct solvers do not allow for a linear time complexity.
	 Nevertheless, they still perform well in some situations \cite{george_complexity_1988}.
	 On the contrary, iterative methods may archive a linear time complexity if the number of iterations is uniformly
	 bounded and each iteration is of linear complexity.
	 Work estimates \cite{bramble_multigrid_1993} prove this for the multigrid $V(r)$-cycle under the assumption of an (asymptotically) exponential growth in the
	 degrees of freedom, i.e., $N_{\ell+1}\geq a N_{\ell}$ for some $a>1$.

	 The numerical benchmarks indeed verify a linear $\mathcal O(N)$ time complexity of the iterative MG and PCG solvers
	 and a superlinear growth of approximately $\mathcal O(N^{1.4})$ for the direct solver.
	 Figure \ref{fig:Timing} verifies the uniform bound on the number of iterations for MG and PCG with the $V(r)$-cycle
	 for $r=1,2,5$.
	 This results in observed optimal convergence rates of the total error $\trb{u-\tilde u_h}$ as well as the error estimator $\eta(\T)$ with
	 respect to the computational time in seconds for the iterative schemes and the experimental convergence rate $1.4$
	 for AFEM with direct solve.
	 For the shown benchmark B3, the PCG method with a single $r=1$ smoothing iteration performs the best and obtains the
	 solution to \eqref{eqn:ADWP} faster than
	 the direct solver on fine meshes with more than $3\times10^4$ degrees of freedom.

\subsection{Norm of the multigrid iteration matrix}
\label{sub:Norm_iteration_matrix}
\begin{figure}[]
		\centering
		\hspace*{-3em}\includegraphics{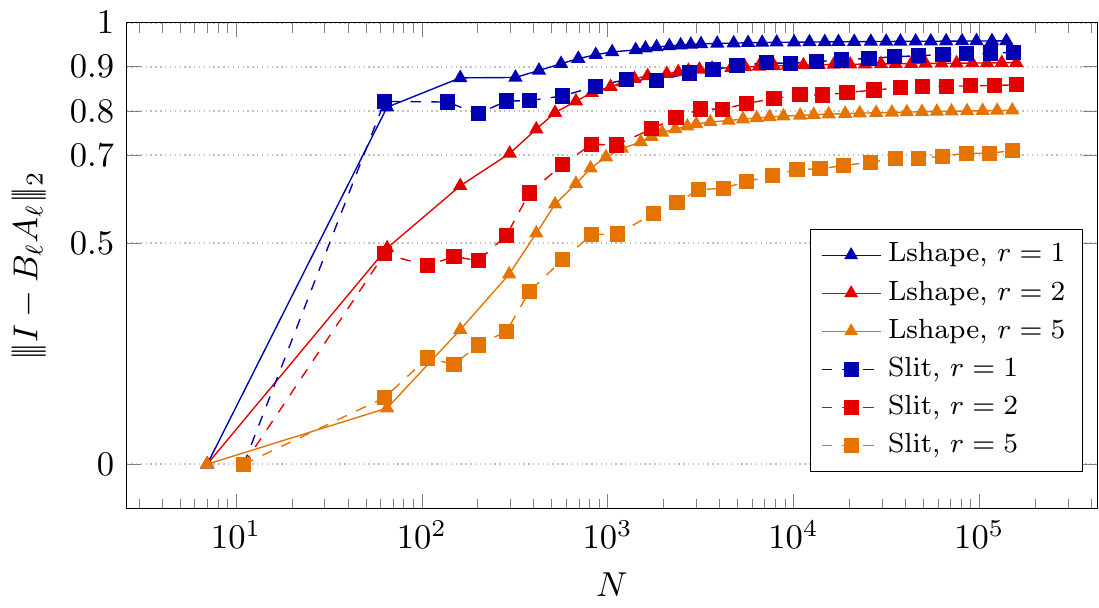}
		\caption{History of $\trb{\opIBAl}_2$ from the multigrid AFEM for B2 (L-shape) and B3 (Slit) from
		subsections \ref{sub:Benchmarks with homogeneous boundary conditions} and \ref{sub:A benchmark with
	inhomogeneous BC}}
		\label{fig:Reduction_factor}

	\end{figure}

The norm of the multigrid iteration matrix $\opIBAl$ quantifies the convergence rate of the multigrid method and control
over it as in theorem \ref{thm:uniform_approx} leads to uniform convergence.
Given $\ell\in\N$, it is known that \begin{align}\label{eqn:Cc}C\coloneqq \trb{\opIBAl}_2 = c/(1+c)\end{align} for some
$c\in \R$, see the proof of theorem \ref{thm:uniform_approx}. %
The Rayleigh-Ritz principle (also known as the min-max principle) for the symmetric matrices $\opAl$ and $\opAl(\opIBAl)$
shows that the spectral energy norm $\trb{\opIBAl}_2$ of the iteration matrix is equal to the maximal eigenvalue of $\opIBAl$.
Figure \ref{fig:Reduction_factor} displays the history of $\trb{\opIBAl}_2$, computed with \texttt{eigs} from the MATLAB standard
library (and $\opBl$ provided as a function handle) and shows that $\trb{\opIBAl}_2$ is clearly bounded away from one.
Since the work in one $V(r)$-cycle is roughly equal to that of $r$ $V(1)$-cycles, it is expected that the number of iterations
correlates anti-proportionally to the number of smoothing steps $r$.
Table \ref{tab:opIBAl} collects the values of $C, c$ from \eqref{eqn:Cc}, and the number of iterations \texttt{n\_it} on a fine mesh from the
multigrid $\Aext$-AFEM with $r$ smoothing steps and verifies this expectation.

\begin{table}[]
	\centering
	\begin{tabular}{l|llc|llc|llc|}
		$r$ & \multicolumn{3}{c|}{Square (B1)} & \multicolumn{3}{c|}{L-shape (B2)} & \multicolumn{3}{c|}{Slit (B3)}\\
			& $C$	& $c$ & \texttt{n\_it} &$C$	& $c$ & \texttt{n\_it}&$C$	& $c$ & \texttt{n\_it}\\\hline
		$1$ & 0.9014	& 9.14 & 2 &0.9590&23.39&12&0.9339&14.13& 8\\
		$2$ & 0.5774& 1.36 & 1 &0.9096&10.06&6&0.8608&6.18&4\\
		$3$ & 0.3701	& 0.58 & 1 &0.8666&6.49&4&0.8070&4.18&3	\\
		$5$ & 0.1949& 0.24 & 1 &0.8022&4.06&3&0.7157&2.52&2	\\
	\end{tabular}
	\caption{Values of $C\coloneqq \trb{\opIBAl}_2$ for the multigrid AFEM ($\theta=0.5$) on meshes with more than $2\cdot 10^5$ degrees of freedom}
	\label{tab:opIBAl}
\end{table} %

	\FloatBarrier
\section{Concluding remarks}%
\label{sub:Conclusion}
	A comparison between the classical and the extended Argyris space in section \ref{sec:Numerical results} shows
	qualitatively similar mesh sequences and convergence.
	The extended Argyris space comes with about $11$\% more degrees of freedom compared to the standard Argyris space on
	the same mesh.
	This extra amount of computation stands opposed to the availability of theoretically justified, local multilevel
	preconditioned solvers and an easily computable, reliable, and efficient estimator of the algebraic error.
Moreover, section \ref{sec:Numerical2} finds the $\Aext$-AFEM algorithm with inexact solve highly efficient and of
linear time complexity in the overall computational cost.
The linear space (memory) complexity of the local multilevel scheme is another advantage over classical direct solvers
\cite{george_complexity_1988}
that becomes important on systems with limited memory.
	A possible extension to the standard Argyris FEM requires a different prolongation as no natural injection from
	coarse to fine spaces exists.
	Theoretical results for the related adaptive algorithm are not available (see \cite{bramble_multigrid_1995} for
	multigrid on quasi-uniform meshes).

	The local Gauß-Seidel smoother only acts on the refined portion of the mesh.
	A simple MATLAB implementation was found competitive with the highly optimised direct solver and the PCG solver is
	already faster on meshes with $3\times10^4$ degrees of freedom.
	The alternative application of the standard Gauß-Seidel (also known as multiplicative) smoother that acts on the full
	set of degrees of freedom shows no qualitative reduction in the number of iterations and comes with an additional
	computational cost.
	This extra cost can be efficiently circumvented with local smoothing, see also \cite{wu_uniform_2006}, considered in
	this paper.

	A possible (uniform) convergence rate of the multigrid iteration close to one does not spoil the
	application of multilevel preconditioned methods in the adaptive setting, as the moderate number of iterations in
	table \ref{tab:opIBAl} suggests.
	In fact, the PCG method with a single smoothing step only required between $1$ and $4$ iterations throughout.

	The hierarchical Argyris FEM marks a paradigm shift in the approximation of conforming fourth-order problems away
	from minimising the dimension of the discrete spaces towards justifying higher order methods.
	Numerical benchmarks reestablish the Argyris element with high convergence rates, an easy
	implementation by transformation \cite{dominguez_algorithm_2008,kirby_general_2017}, and an overall linear time
	complexity of the optimal adaptive algorithm.

\FloatBarrier
\section*{Acknowledgements}
This work was supported by the
DFG Priority Program~1748
\emph{Reliable Simulation Techniques in Solid Mechanics.\ Development of
Non-standard Discretization Methods, Mechanical and Mathematical Analysis}
within the project
\emph{Foundation and application of generalized mixed FEM towards
nonlinear problems in solid mechanics} (CA 151/22-2) and  the \emph{Berlin Mathematical School}.
Furthermore, the author gratefully acknowledges the valuable advice of Prof.\ Carsten Carstensen from the Humboldt
Universität zu Berlin as well as his supervision throughout the studies.
		\bibliographystyle{elsarticle-harv}
		\bibliography{./Bibliography}
\end{document}